\documentclass[12pt,a4paper]{article}
\usepackage{amssymb}
\usepackage{amsfonts}

\usepackage{amsmath}
\usepackage{amsthm}
\usepackage{enumerate}

\RequirePackage[OT1]{fontenc}
\RequirePackage{amsthm,amsmath, xcolor}
\RequirePackage[numbers]{natbib}
\RequirePackage[colorlinks,citecolor=blue,urlcolor=blue]{hyperref}
\usepackage{amscd,amsfonts,amssymb,latexsym,array,hhline,graphicx}
\usepackage{float}

\newtheorem{theorem}{Theorem}[section]
\newtheorem{proposition}[theorem]{Proposition}

\newtheorem{lemma}[theorem]{Lemma}

\newtheorem{remark}{Remark}[section]
\newtheorem{corollary}[theorem]{Corollary}

\newcommand\cA{{\cal A}}

\newcommand\cG{{\cal G}}
\newcommand\cF{{\cal F}}

\newcommand\cB{{\cal B}}

\newcommand\cN{{\cal N}}

\newcommand\cX{{\cal X}}
\newcommand\cD{{\cal D}}
\newcommand\cQ{{\cal Q}}

\newcommand\cR{{\cal R}}

\newcommand\ve{\varepsilon}

\newcommand\Er{\mbox{Err}}

\def\bbr{{\mathbb R}}

\def\text#1{\hbox{#1}}
\def\proof{{\noindent \bf Proof. }}
\def\endproof{\mbox{\ $\qed$}}

\def\E{{\bf E}}

\def\P{{\bf P}}
\def\B{{\bf B}}
\def\p{{\bf p}}
\def\g{{\bf g}}
\def\l{{\bf l}}
\def\C{{\bf C}}
\def\D{{\bf D}}

\def\G{{\bf G}}
\def\L{{\bf L}}
\def\U{{\bf U}}
\def\M{{\bf M}}

\def\l{{\bf l}}

\def\R{{\bf R}}

\def\Chi{{\bf 1}}

\def\d{\mathrm{d}}
\def\build #1_#2{\mathrel{\mathop{\kern 0pt #1}\limits_{#2}}}
\newcommand\tr{\mbox{tr}}
\newcommand\Trg{\mbox{Tr}}
\newcommand{\wh}{\widehat}
\newcommand{\wt}{\widetilde}
\newcommand{\zs}[1]{{\mathchoice{#1}{#1}{\lower.25ex\hbox{$\scriptstyle#1$}}
{\lower0.25ex\hbox{$\scriptscriptstyle#1$}}}}

\numberwithin{equation}{section}

\begin{document}
\title{
Improved
robust
model selection
methods
for a L\'evy
nonparametric
 regression
 in continuous time
\thanks{
This work was supported by
 the RSF grant 17-11-01049 (National Research Tomsk State University).
}}
\author{Pchelintsev E.A.,
\thanks{
Department of Mathematical Analysis and Theory of Functions,
 Tomsk State University,
Lenin str. 36,
 634050 Tomsk, Russia,
 e-mail: evgen-pch@yandex.ru
}
\and
Pchelintsev V.A.,
\thanks{
Department of Mathematics and Informatics,
 Tomsk Polytechnic University,
Lenin str. 30,
 634050 Tomsk, Russia,
 e-mail: vpchelintsev@vtomske.ru
}
 \and
 Pergamenshchikov S.M.\thanks{
 Laboratoire de Math\'ematiques Raphael Salem,
 Avenue de l'Universit\'e, BP. 12,
  Universit\'e de Rouen,
   F76801, Saint Etienne du Rouvray, Cedex France
   and
International Laboratory of Statistics of Stochastic Processes and
Quantitative Finance of Tomsk State University,
 e-mail: Serge.Pergamenchtchikov@univ-rouen.fr}
}
 \date{}
\maketitle

\begin{abstract}
In this paper we
develop the James -  Stein improved  method for
the estimation problem of
 a nonparametric
periodic function observed
 with  L\'evy noises in continuous time.
An adaptive model selection
procedure based on the  weighted improved least squares estimates is constructed.
The improvement effect  for nonparametric models is studied.
It turns out that in  non-asymptotic setting the accuracy improvement for nonparametric models is more important,
 than for parametric ones.
Moreover,  sharp oracle inequalities for the robust
risks have been shown and the adaptive efficiency property for the proposed procedures
has been established. The numerical simulations are given.
\end{abstract}

\\
\par
{\bf AMS (2010) Subject Classification : primary 62G08; secondary 62G05}

\bibliographystyle{plain}


\section{Introduction}\label{sec:In}

Consider the following
nonparametric regression model in continuous time
\begin{equation}\label{sec:In.1}
 \d\,y_\zs{t} = S(t)\d\,t + \,
\d \xi_\zs{t}\,,\quad
 0\le t \le n\,,
\end{equation}
where $S(\cdot)$ is an unknown $1$ - periodic
 function,
$(\xi_\zs{t})_\zs{0\le t\le n}$ is an unobserved noise. The problem is to estimate
the function $S$ on the observations $(y_\zs{t})_\zs{0\le t\le n}$.
 Note that, if $(\xi_\zs{t})_\zs{0\le t\le n}$ is a  Brownian motion, then we obtain the well-known "signal+white noise" model
which is very popular in statistical radio-physics (see, for example, \cite{IbragimovKhasminskii1981, Kutoyants1977, Kutoyants1984, Pinsker1981}).
 In this paper we assume that in addition
 to intrinsic noises in  radio-electronic systems, approximated usually by the gaussian white or color noise,
 the useful signal $S$ is distorted by the impulse flow described
by  L\'evy processes
defined in the next section.
The cause
of a pulse stream
 can be, for example, either external unintended (atmospheric) or intentional impulse noises
 or
 errors in the demodulation and the channel decoding for  binary information symbols.
 Note that, for the first time the impulse noises
for the detection signal problems
 have been studied
 by Kassam in
\cite{Kassam1988}
 through  compound Poisson processes.
 Later, such processes was used
  in
  \cite{Flaksman2002, KonevPergamenshchikov2012, KonevPergamenshchikov2015, KPP2014, Pchelintsev2013}
for
parametric and
nonparametric signal estimation problems.
It should be noted that such models
are too limited, since   the compound Poisson process can  describe  only the  large impulses influence  with a single fixed  frequency.
However, the real technical (for example, telecommunication or navigation) systems
work under  noise impulses having different sizes and different frequencies (see, for example, \cite{Proakis1995}).
To take this into account
 one needs to use many (may be infinite number)  different compound Poisson processes in the same observation model.
This is possible to do only in a framework of  L\'evy processes which are natural extensions for  the
 compound Poisson processes. Moreover, it should be
noted also that
L\'evy  models are
fruitfully
used in the different applied problems
(see, for example,
\cite{BarndorffNielsenShephard2001, Bertoin1996, ComteGenenCatalot2011, ContTankov2004} and the references therein).
In this paper we consider the adaptive estimation problem for the function $S$
 i.e. when its regularity properties are unknown.  To do this  we use the model selection methods.
The interest to such statistical procedures is explained by the fact that
they provide
adaptive solutions for the nonparametric estimation through
oracle inequalities which give  non-asymptotic upper
bounds for the quadratic risks including
the minimal risk over chosen  the estimators family.
It will be noted that for the first time the model selection methods
were proposed by
 Akaike \cite{Akaike1974} and Mallows \cite{Mallows1973}
for parametric models. Then,
these methods have been developed
for  nonparametric estimation problems
 by
Barron, Birg\'e and Massart \cite{BarronBirgeMassart1999}
and
Fourdrinier and Pergamenshchikov \cite{FourdrinierPergamenshchikov2007}
for regression models in discrete time
and
Konev and Pergamenshchikov
\cite{KonevPergamenshchikov2010} in continuous time.
Unfortunately, the oracle inequalities obtained in these papers
 can not
provide the efficient estimation in the adaptive setting, since
the upper bounds in these inequalities
have some fixed coefficients in the main terms which are more than one.
To obtain  the efficiency property
 one has to obtain
 the sharp oracle inequalities, i.e. the inequalities
 in which
  the coefficient
   at the principal term
 is close to unity. To obtain such inequalities
 for  general non-Gaussian observations one needs to use the
 method
 proposed by
Konev and Pergamenshchikov in
\cite{KonevPergamenshchikov2009a, KonevPergamenshchikov2009b, KonevPergamenshchikov2012, KonevPergamenshchikov2015}
for semimartingale models in  continuous time based on the
model selection tool
developed
  by
 Galtchouk and Pergamenshchikov in
 \cite{GaltchoukPergamenshchikov2009a, GaltchoukPergamenshchikov2009b}  for  heteroscedastic
 non-Gaussian
 regression models in discrete time.

The goal of this paper is to develop a new sharp model selection method for estimating the unknown signal $S$
using the improved estimation approach.
Usually, the model selection procedures are based on the least squares estimates.  This paper proposes
 the improved least squares estimates which enable us to improve considerably the non-asymptotic estimation accuracy.
For the first time such idea  was proposed
by Fourdrinier and Pergamenshchikov
 in \cite{FourdrinierPergamenshchikov2007}
for regression models in discrete time and
by Konev and Pergamenshchikov
in
\cite{KonevPergamenshchikov2010}
for  Gaussian regression models in continuous time. We develop these methods for the non-Gaussian
regression models
in continuous time.
It should be noted that generally for the conditionally Gaussian regression models  we can not use the well-known improved estimators proposed in
\cite{FourdrinierStrawderman1996, JamesStein1961} for Gaussian or spherically symmetric observations. To apply the improved estimation methods to the
 non-Gaussian regression models in continuous time
 one needs to use the modifications of the well-known James - Stein estimators  proposed  in  \cite{KPP2014, Pchelintsev2013} for parametric problems.
We use these estimators to construct  model selection procedures for nonparametric models. Then to study the non-asymptotic accuracy
we develop a special analytical tool for the L\'evy regression models   to obtain sharp oracle inequalities
for the improved model selection procedures.
Then to study the efficiency property for the proposed estimation procedure we need to obtain a lower bound
for the quadratic risks. Usually, to do this one uses the van Trees inequality.  In this paper we show
the corresponding van Trees inequality
for the L\'evy  regression models  and then
 we derive the needed asymptotic sharp lower bound for the normalized risks, i.e. we find the Pinsker constant
 for the model   \eqref{sec:In.1}.
  As to the upper bound,
 similarly to
\cite{KonevPergamenshchikov2009b},
  we use the obtained  sharp oracle inequality for the
 weighted least squares estimators  containing the efficient Pinsker procedure. Therefore, through the oracle inequality
 we estimate from above the risk of the proposed  procedure by the risk
of the efficient Pinsker procedure up to some coefficient which goes to one. As a result we show the asymptotic efficiency
without using the smoothness information of the function $S$.

The rest of the paper is organized as follows.
In Section \ref{sec:Mod} we describe the noise processes in  \eqref{sec:In.1} and define the main
risks for the estimation problem.
In Section  \ref{sec:Imp} we construct the improved least squares estimates and
study the improvement effect for the L\'evy model. In Section \ref{sec:Mo} we construct the improved model selection
procedure and show the sharp oracle inequalities. In Section \ref{sec:Sim} the Monte Carlo simulation results are given.
 The asymptotic efficiency is studied in Section \ref{sec:Ae}. In Section \ref{sec:Stc}
we study some properties of the stochastic integrals with respect to the L\'evy processes.
In Section \ref{sec:VanTrees} we prove the van Trees inequality for the model \eqref{sec:In.1}.
 In Section \ref{sec:Prf} we prove all main results and in Appendix we give some technical results.

\bigskip

\section{Noise process model}\label{sec:Mod}

 First, we assume that
 the noise process $(\xi_\zs{t})_\zs{0\le t\le n}$ in \eqref{sec:In.1}
is defined as
\begin{equation}\label{sec:In.1+1}
\xi_\zs{t}=\sigma_\zs{1} w_\zs{t} + \sigma_\zs{2} z_\zs{t}
\quad\mbox{and}\quad
z_\zs{t}=x*(\mu-\wt{\mu})_\zs{t}
\,,
\end{equation}
where $\sigma_\zs{1}$ and $\sigma_\zs{2}$ are  some unknown constants,
$(w_\zs{t})_\zs{t\ge\,0}$ is a standard Brownian motion,
 "$*$"\ denotes the stochastic integral with respect to the compensated jump measure $\mu(\d s\,\d x)$ with  deterministic
compensator $\wt{\mu}(\d s\,\d x)=\d s\Pi(\d x)$, i.e.
$$
z_\zs{t}=\int_0^t\int_{\bbr_\zs{0}}x\,(\mu-\wt{\mu})(\d s \,\d x)\,.
$$
Here $\Pi(\cdot)$ is a L\'evy measure, i.e.  some positive measure on $\bbr_\zs{0}=\bbr\setminus \{0\}$,
(see, for example,
\cite{ContTankov2004, JacodShiryaev2002} for details) such that
\begin{equation}\label{sec:Ex.1-00_mPi}
\Pi(x^{2})=1
\quad\mbox{and}\quad
\Pi(x^{6})
\,<\,\infty\,.
\end{equation}
We use the notation $\Pi(\vert x\vert^{m})=\int_\zs{\bbr_\zs{0}}\,\vert z\vert^{m}\,\Pi(\d z)$. Note that the L\'evy measure
 $\Pi(\bbr_\zs{0})$ may be equal to $+\infty$. It should be noted that in all papers on
 the nonparametric signal estimation in the model \eqref{sec:In.1}
  the main condition on the jumps  is the finiteness of the L\'evy measure, i.e.  $\Pi(\bbr_\zs{0})<+\infty$.

The process \eqref{sec:In.1+1} allows us to consider the
several independent  impulse noise sources with the different  frequencies.
Indeed, in this case (see, for example, page 135 in \cite{ContTankov2004}) we introduce compound Poisson processes into the model \eqref{sec:In.1}  as
\begin{equation}\label{sec:In.1_ex_1}
z_\zs{t}=\sum^{M}_\zs{k=1}\sum^{N^{k}_\zs{t}}_\zs{j=1}\,Y_\zs{k,j}
\,,
\end{equation}
where $(N^{1}_\zs{t})_\zs{t\ge 0},\ldots,(N^{M}_\zs{t})_\zs{t\ge 0}$ are independent Poisson processes
with the intensities $\lambda_\zs{1},\ldots, \lambda_\zs{M}$
and the sizes of impulses $(Y_\zs{1,j})_\zs{j\ge 1},\ldots, (Y_\zs{M,j})_\zs{j\ge 1}$ are
independent i.i.d. sequences with $\E Y_\zs{k,j}=0$ and $\varsigma^{2}_\zs{k}=\E Y^{2}_\zs{k,j}<\infty$.
In this case the L\'evy measure for any Borel set $\Gamma\subseteq\bbr_\zs{0}$ is defined as
$$
\Pi(\Gamma)=\sum^{M}_\zs{k=1}\lambda_\zs{k}\,\P(Y_\zs{k,1}\in\Gamma)\,.
$$
Next, note, that if
\begin{equation}\label{sec:In.1_cond_Ex}
\sum_\zs{k\ge 1}\,\lambda_\zs{k}\,\varsigma^{2}_\zs{k}
<\infty\,,
\end{equation}
then we can introduce the infinite number of the noise jumps setting
\begin{equation}\label{sec:In.1_ex_2}
z_\zs{t}=\sum^{\infty}_\zs{k=1}\sum^{N^{k}_\zs{t}}_\zs{j=1}\,Y_\zs{k,j}
\,.
\end{equation}
Moreover, if the total noise intensity $\sum_\zs{k\ge 1}\lambda_\zs{k}=+\infty$, then  $\Pi(\bbr_\zs{0})=+\infty$,
i.e. we obtain the observation model  with saturated impulse noise.

In the sequel we will denote by $Q$ the distribution of the process $(\xi_\zs{t})_\zs{0\le t\le n}$
in the Skorokhod space $\D[0,n]$ and by $\cQ_\zs{n}$ we denote all these distributions for which
  the parameters
    $\sigma_\zs{1}$
and $\sigma_\zs{2}$ satisfy the conditions
\begin{equation}\label{sec:Ex.01-1}
0< \underline{\sigma}\le \sigma^{2}_\zs{1}
\quad\mbox{and}\quad
\sigma=\sigma^{2}_\zs{1}+\sigma^{2}_\zs{2}\,
\le
\overline{\sigma}
\,,
\end{equation}
where
 the bounds
 $\underline{\sigma}$
and
$\overline{\sigma}$
are functions of $n$, i.e.
 $\underline{\sigma}= \underline{\sigma}_\zs{n}$
and
$\overline{\sigma}=\overline{\sigma}_\zs{n}$
 such that for any $\epsilon>0$
\begin{equation}\label{sec:Ex.01-2}
\liminf_\zs{n\to\infty}\,n^{\epsilon}\,
\underline{\sigma}_\zs{n}
\,
>0
\quad\mbox{and}\quad
\lim_\zs{n\to\infty}\,n^{-\epsilon}\,\overline{\sigma}_\zs{n}
=0
\,.
\end{equation}

We also assume that the
distribution $Q$
of
 the noise process  $(\xi_\zs{t})_\zs{0\le t\le n}$ is  unknown.
 We know only that this distribution belongs to the
  distribution family $\cQ_\zs{n}$ defined  in
  \eqref{sec:Ex.01-1}--\eqref{sec:Ex.01-2}.
    By these reasons  we use   the robust estimation approach developed for nonparametric problems in
 \cite{GaltchoukPergamenshchikov2006, KonevPergamenshchikov2012, KonevPergamenshchikov2015}.
To this end we will measure the estimation quality by the robust risk defined as
\begin{equation}\label{sec:risks_11_0}
\cR^{*}_\zs{n}(\wh{S}_\zs{n},S)=\sup_\zs{Q\in\cQ_\zs{n}}\,
\cR_\zs{Q}(\wh{S}_\zs{n},S)\,,
\end{equation}
where  $\wh{S}_\zs{n}$ is an estimate, i.e. any function of $(y_\zs{t})_\zs{0\le t\le n}$,
$\cR_\zs{Q}(\cdot,\cdot)$ is the usual quadratic risk defined as
\begin{equation}\label{sec:risks_00}
\cR_\zs{Q}(\wh{S}_\zs{n},S):=
\E_\zs{Q,S}\,\|\wh{S}_\zs{n}-S\|^2
\quad\mbox{and}\quad
\Vert S\Vert^{2}=\int^{1}_\zs{0}\,S^{2}(t)\d t
\,.
\end{equation}
The first goal in this paper is to
develop shrinkage nonparametric estimation methods for  $S$
which improve the non asymptotic robust estimation accuracy \eqref{sec:risks_11_0}
with respect to the well known least squares estimators.
The next goal is to provide non asymptotic optimality in the sense of sharp oracle inequalities.
 Moreover, asymptotically, as $n\to\infty$,  our goal is to show the efficiency property
for the proposed shrinkage estimators for the risks
\eqref{sec:risks_11_0}.

\bigskip

\section{Improved estimation}\label{sec:Imp}

Let $(\phi_\zs{j})_\zs{j\ge\, 1}$ be an orthonormal basis in $\L_\zs{2}[0,1]$.
We extend these functions  by the periodic way on $\bbr$, i.e.  $\phi_\zs{j}(t)$=$\phi_\zs{j}(t+1)$ for any $t\in\bbr$.

$\B_\zs{1}$) {\em Assume that the basis functions are uniformly bounded, i.e.
for some  $\overline{\phi}_\zs{n}>0$, which in general case may be depend on $n$,
}
\begin{equation}\label{sec:In.3-00_Upb}
\sup_\zs{0\le j\le n}\,\sup_\zs{0\le t\le 1}\vert\phi_\zs{j}(t)\vert\,
\le\,
\overline{\phi}_\zs{n}
<\infty\,.
\end{equation}

\noindent  For example, we can take
 the trigonometric basis    defined as $\Trg_\zs{1}\equiv 1$ and for $j\ge 2$
\begin{equation}\label{sec:In.5_Trb}
 \Trg_\zs{j}(x)= \sqrt 2
\left\{
\begin{array}{c}
\cos(2\pi[j/2]x)\, \quad\mbox{for even}\quad j \,;\\[4mm]
\sin(2\pi[j/2]x)\quad\mbox{for odd}\quad j\,,
\end{array}
\right.
\end{equation}
where  $[a]$ denotes integer part of $a$.

\bigskip

For estimating the unknown function $S$ in \eqref{sec:In.1} we consider it's Fourier expansion
$$
S(t)=\sum_{j=1}^\infty \theta_\zs{j}\phi_j(t).
$$
The corresponding Fourier coefficients
\begin{equation*}\label{sec:Imp.2}
\theta_\zs{j}=(S,\phi_j)= \int^1_\zs{0}\,S(t)\,\phi_\zs{j}(t)\,\d t
\end{equation*}
can be estimated as
\begin{equation*}\label{sec:Imp.3}
\wh{\theta}_\zs{j}= \frac{1}{n}\int^n_\zs{0}\,\phi_j(t)\,\d
y_\zs{t}\,.
\end{equation*}
In view of \eqref{sec:In.1}, we obtain
\begin{equation}\label{sec:Imp.4}
\wh{\theta}_\zs{j}=\theta_\zs{j}+\frac{1}{\sqrt{n}}\xi_\zs{j}\,,
\end{equation}
where
$$
 \xi_\zs{j}=\frac{1}{\sqrt{n}}
I_\zs{n}(\phi_\zs{j})
\quad\mbox{and}\quad
I_\zs{n}(f)=\int^n_\zs{0}\,f(t)\,\d \xi_\zs{t}\,.
$$
As in \cite{KonevPergamenshchikov2009b} we define a class of weighted least squares estimates for $S(t)$
\begin{equation}\label{sec:Imp.5}
\wh{S}_\zs{\lambda}=\sum^{n}_\zs{j=1}\lambda(j)\wh{\theta}_\zs{j}\phi_\zs{j}\,,
\end{equation}
where the weights $\lambda=(\lambda(j))_{1\leq j\leq n}\in\bbr^{n}$ belong to some finite set $\Lambda$ from $[0,\,1]^n$ for which we set
\begin{equation}\label{sec:Imp.6+0}
\nu_n=\mbox{card}(\Lambda)
\quad\mbox{and}\quad
\vert \Lambda\vert_\zs{n}=\max_{\lambda\in \Lambda} \,L(\lambda)\,,
\end{equation}
where $\mbox{card}(\Lambda)$ is the number of the vectors $\lambda$ in $\Lambda$ and $L(\lambda)=\sum_{j=1}^n\lambda(j)$.
In the sequel we assume that all vectors from $\Lambda$ satisfies the following condition.

\bigskip

$\B_\zs{2})$ {\sl Assume that
for any vector $\lambda\in\Lambda$ there exists
 some fixed integer $d=d(\lambda)$
 such that their first $d$ components
equal to one, i.e. $\lambda(j)=1$ for $1\le j\le d$ for any $\lambda\in\Lambda$. }

\begin{remark}\label{Re;sec:Imp.1++}
Note that the weight coefficients satisfying the
  condition $\B_\zs{2})$ was introduced in
\cite{Nussbaum1985} to
construct the efficient estimation
for the nonparametric
regression model in discrete time.
\end{remark}

Now we need the $\sigma$ - field generated by the jumps of the process \eqref{sec:In.1+1}, i.e. we set
$\cG_\zs{n}=\sigma\{z_\zs{t}\,,0\le t\le n\}$. To construct the improved estimators we need the following proposition.

\begin{proposition}\label{sec:Imp.Prop_2_1}
For any $n\ge 1$
 the random vector
 $\wt{\xi}_\zs{d,n}=(\xi_\zs{j})_\zs{1\le j\le d}$
 is the $\cG_\zs{n}$ - conditionally Gaussian in $\bbr^{d}$ with zero mean and the covariance matrix
\begin{equation*}\label{sec:Imp.6-1}
\G_\zs{n}=\left(
\E\,\xi_\zs{i}\,\xi_\zs{j}|\cG_\zs{n}
\right)_\zs{1\le i,j\le d}
\end{equation*}
such that
\begin{equation}\label{sec:Imp.6-1-0}
\inf_\zs{Q\in\cQ_\zs{n}}
\left(\tr \,\G_\zs{n}-\lambda_\zs{max}(\G_\zs{n})\right)
\geq (d-1)\underline{\sigma}_\zs{n}\,,
\end{equation}
where $\lambda_\zs{max}(A)$ is the maximal eigenvalue of the matrix $A$.
\end{proposition}

Now,
for the first $d$ Fourier coefficients in \eqref{sec:Imp.4}
 we use
 the improved estimation method proposed
for parametric models
 in
\cite{Pchelintsev2013}. To this end we
set $\wt{\theta}_\zs{n}=(\wh{\theta}_\zs{j})_\zs{1\le j\le d}$.
In the sequel we will use the norm $\vert x\vert^{2}_\zs{d}=\sum^{d}_\zs{j=1}\,x^{2}_\zs{j}$
for any vector $x=(x_\zs{j})_\zs{1\le j\le d}$ from $\bbr^{d}$.
Now
we define the shrinkage estimators as
\begin{equation}\label{sec:Imp.12_Imp}
\theta^{*}_\zs{\lambda,j}=
\left(1-g_\zs{\lambda}(j)\right)\wh{\theta}_\zs{j}
\quad\mbox{and}\quad
g_\zs{\lambda}(j)=
\frac{c_\zs{n}}{|\wt{\theta}_\zs{n}|_\zs{d}}
 \Chi_\zs{\{1\le j\le d\}}
\,,
\end{equation}
where
$$
c_\zs{n}=
c_\zs{n}(\lambda)
=
\frac{(d-1)\underline{\sigma}_\zs{n}}{\left(r_\zs{n}+\sqrt{d\,\overline{\sigma}_\zs{n}/n}\right)\,n}
$$
and the threshold $\overline{\sigma}_\zs{n}>0$ is given in the lower bound \eqref{sec:Imp.6-1-0}.
The positive parameter
 $r_\zs{n}$  is such that
\begin{equation}\label{sec:Imp.12+1_r_n}
\lim_\zs{n\to\infty}\,r_\zs{n}\,=\infty
\quad\mbox{and}\quad
\lim_\zs{n\to\infty}\,
\frac{r_\zs{n}}{n^{\epsilon}}
\,=\,0
\end{equation}
for any $\epsilon>0$.

Now we introduce a class of shrinkage
weighted least squares estimates for $S$ as
\begin{equation}\label{sec:Imp.11}
S^{*}_\zs{\lambda}=\sum^{n}_\zs{j=1}\lambda(j)\theta^{*}_\zs{\lambda,j}\phi_\zs{j}\,.
\end{equation}

We denote the difference of quadratic risks of the estimates \eqref{sec:Imp.5} and \eqref{sec:Imp.11} as
$$
\Delta_{Q}(S):=\cR_\zs{Q}(S^{*}_\zs{\lambda},S)-\cR_\zs{Q}(\wh{S}_\zs{\lambda},S)\,.
$$
We obtain the following result.

\begin{theorem}\label{Th.sec:Imp.1}
Let the observed process $(y_\zs{t})_\zs{0\le t \le n}$ describes by the equations \eqref{sec:In.1}--\eqref{sec:In.1+1}. Then for any $n\ge 1$
\begin{equation}\label{sec:Imp.11+1}
\sup_{Q\in\cQ_\zs{n}}\,\sup_\zs{\Vert S\Vert\le r_\zs{n}}
\Delta_{Q}(S)\le-c^2_\zs{n}
\,.
\end{equation}
\end{theorem}

\begin{remark}\label{Re;sec:Imp.1}
The inequality \eqref{sec:Imp.11+1} means that non-asymptotically, i.e. for non large $n\ge 1$,
 the estimate \eqref{sec:Imp.11}   outperforms in mean square accuracy the estimate \eqref{sec:Imp.5}.
As we will see later in the efficient weight coefficients  $d \approx n^{\epsilon}$ as $n\to \infty$ for some $\epsilon>0$.
Therefore, in view of the  definition of the constant $c_\zs{n}$
in
\eqref{sec:Imp.12_Imp}
 and the conditions
\eqref{sec:Ex.01-2}
 and
 \eqref{sec:Imp.12+1_r_n}
 $n c_\zs{n}\to\infty$ as $n\to\infty$.
This means that improvement is considerably may better than for the parametric regression
when the parameter dimension $d$ is fixed  \cite{Pchelintsev2013}.
\end{remark}

\bigskip

\section{Model selection}\label{sec:Mo}

In this section we construct  a model selection procedure  for
the estimation of  $S$ in \eqref{sec:In.1} on the basis of the weighted shrinkage estimators \eqref{sec:Imp.11}.
To this end we consider  the empirical squared error defined as
$$
\Er_\zs{n}(\lambda)=\|S^*_\zs{\lambda}-S\|^2.
$$
In order to obtain a good estimate, we have to write a rule to choose a weight vector
$\lambda\in \Lambda$ in \eqref{sec:Imp.11}. It is obvious, that the best way is to minimise
the empirical squared error with respect to $\lambda$. Making use the estimate definition
\eqref{sec:Imp.11} and the Fourier transformation of $S$ implies
\begin{equation}\label{sec:Mo.1}
\Er_\zs{n}(\lambda)\,=\,
\sum^{n}_\zs{j=1}\,\lambda^2(j)(\theta^{*}_\zs{\lambda,j})^2\,-
2\,\sum^{n}_\zs{j=1}\,\lambda(j)\theta^{*}_\zs{\lambda,j}\,\theta_\zs{j}
+
\Vert S\Vert^{2}\,.
\end{equation}
Since the Fourier coefficients $(\theta_\zs{j})_\zs{j\ge 1}$ are
unknown, the weight coefficients $(\lambda_\zs{j})_\zs{j\ge 1}$ can
not be found by minimizing this quantity. To circumvent this
difficulty one needs to replace  the terms
$\theta^{*}_\zs{\lambda,j}\,\theta_\zs{j}$ by their estimators
$\overline{\vartheta}_\zs{\lambda,j}$ defined as
\begin{equation}\label{sec:Mo.2vartheta}
\overline{\vartheta}_\zs{\lambda,j}=
\theta^{*}_\zs{\lambda,j}\,\wh{\theta}_\zs{j}-\frac{\wh{\sigma}_\zs{n}}{n}\,,
\end{equation}
where $\wh{\sigma}_\zs{n}$ is the estimate for the limiting variance
of $\sigma=\E_\zs{Q}\,\xi^{2}_\zs{j}$ which we choose in the following form
\begin{equation}\label{sec:Mo.3}
\wh{\sigma}_\zs{n}=\sum_\zs{j=[\sqrt{n}]+1}^n \wh{t}_\zs{j}^2
\quad\mbox{and}\quad
\wh{t}_\zs{j}=
\frac{1}{n}\,
\int_0^{n} \Trg_\zs{j}(t)\d y_t
\,.
\end{equation}
For this change in the empirical squared error, one has to pay
some penalty. Thus, one comes to the cost function of the form
\begin{equation}\label{sec:Mo.4}
J_\zs{n}(\lambda)\,=\,\sum^{n}_\zs{j=1}\,\lambda^2(j)(\theta^{*}_\zs{\lambda,j})^2\,-
2\,\sum^{n}_\zs{j=1}\,\lambda(j)\,\overline{\vartheta}_\zs{\lambda,j}
+\,\delta\,\wh{P}_\zs{n}(\lambda)\,,
\end{equation}
where $\delta$ is some positive constant,
$\wh{P}_\zs{n}(\lambda)$ is the penalty term defined as
\begin{equation}\label{sec:Mo.5_whP}
\wh{P}_\zs{n}(\lambda)=\frac{\wh{\sigma}_\zs{n}\,|\lambda|^2_\zs{n}}{n}
\,.
\end{equation}

\noindent
We define the improved model selection procedure as
\begin{equation}\label{sec:Mo.6}
S^*=S^*_\zs{\lambda^*}
\quad\mbox{and}\quad
\lambda^*=\mbox{argmin}_\zs{\lambda\in\Lambda}\,J_n(\lambda)\,.
\end{equation}
It will be noted that $\lambda^*$ exists because
 $\Lambda$ is a finite set. If the
minimizing sequence in \eqref{sec:Mo.6} $\lambda^*$ is not
unique, one can take any minimizer. Now, to write the oracle inequality we set
\begin{equation}\label{sec:Mo.9+1_psi}
\Psi_\zs{Q,n}=(1+\overline{\phi}_\zs{n}^4)\,(1+\sigma)(1+c^{*}_\zs{n})\nu_n\,,
\end{equation}
where  $c^{*}_\zs{n}=n\max_\zs{\lambda\in \Lambda}c^{2}_\zs{n}(\lambda)$.
It is useful to note that in view of the first condition in \eqref{sec:Ex.01-2}
and the properties \eqref{sec:Imp.12+1_r_n}
 the constant $c^{*}_\zs{n}$ is not large as $n\to\infty$, i.e.
 for any $\epsilon>0$
\begin{equation}\label{sec:Mo.8_cn}
\lim_\zs{n\to\infty}\,\frac{c^{*}_\zs{n}}{n^{\epsilon}}
=0\,.
\end{equation}

\noindent First we study the non asymptotic properties for the procedure \eqref{sec:Mo.6}.

\begin{theorem}\label{sec:Mo.Th.1}
There exists some constant $\check{\l}>0$ such that
for any $n\geq1$ and $0<\delta<1/2$, the risk \eqref{sec:risks_00} of estimate \eqref{sec:Mo.6} for $S$
satisfies the oracle inequality
\begin{align}\nonumber
\cR_\zs{Q}(S^*,S)\,&\le\, \frac{1+5\delta}{1-\delta}
\min_\zs{\lambda\in\Lambda} \cR_\zs{Q}(S^*_\zs{\lambda},S)
+\check{\l}\frac{\Psi_\zs{Q,n}}{n \delta}\,
\\[2mm]\label{sec:Mo.10}
&
+
\frac{12\vert\Lambda\vert_\zs{n}\E_\zs{Q}|\wh{\sigma}_\zs{n}-\sigma|}{n}\,.
\end{align}
\end{theorem}

\noindent
 In the case, when the value of $\sigma$ is known, one can take
$\wh{\sigma}_\zs{n}=\sigma$ and
\begin{equation}\label{sec:Mo.9_P}
P_\zs{n}(\lambda)=\frac{\sigma\,|\lambda|^2_\zs{n}}{n}\,,
\end{equation}
then we can rewrite the oracle inequality \eqref{sec:Mo.10}  in the following form
\begin{equation*}\label{sec:Mo.10+1}
\cR_\zs{Q}(S^*,S)\,\le\, \frac{1+5\delta}{1-\delta}
\min_\zs{\lambda\in\Lambda} \cR_\zs{Q}(S^*_\zs{\lambda},S)
+\check{\l}\frac{\Psi_\zs{Q,n}}{n \delta}\,
\,.
\end{equation*}
Also we study the accuracy properties for the estimator \eqref{sec:Mo.3}.

\begin{proposition}\label{sec:Mo.Prop.1}
Let in the model \eqref{sec:In.1} the function $S(\cdot)$ is continuously differentiable.
Then, there exists some constant $\check{\l}>0$ such that
for any $n\geq2$ and $S$
\begin{equation}\label{sec:Mo.15++-00}
\E_\zs{Q}|\wh{\sigma}_\zs{n}-\sigma| \leq
\,\check{\l}\,
\frac{(1+\|\dot{S}\|^2)}{\sqrt{n}}\,,
\end{equation}
where $\dot{S}$ is the derivative of the function $S$.
\end{proposition}

\begin{remark}\label{Re.sec:Ma.TrgBasis}
It should be noted that to estimate the parameter $\sigma$
in \eqref{sec:Mo.2vartheta}
we use the equality \eqref{sec:Imp.4}
 for the Fourier coefficients $t_\zs{j}=(S,\Trg_\zs{j})$ with respect to the trigonometric basis
\eqref{sec:In.5_Trb}, since, as is shown in
Lemma A.6 in \cite{KonevPergamenshchikov2009a}
for any continuously differentiable function $S$
and
for any $m\ge 1$
the sum $\sum_\zs{j\ge m}\,t^{2}_\zs{j}$ can be estimated from above in an explicite form.
Therefore, through the trigonometric basis we can estimate the variance $\sigma$ uniformly over
the functions $S$, when we will study the efficiency property for the proposed procedures.
\end{remark}

\noindent
To obtain the oracle inequality for the robust risk
 we impose the following additional conditions.

\noindent $\C_\zs{1})$ {\it Assume that
 the upper bound for the basic function defined
in \eqref{sec:In.3-00_Upb} is such that for any $\epsilon>0$
 \begin{equation*}\label{sec:Mo.8+2-0}
\lim_\zs{n\to\infty}\,\frac{\overline{\phi}_\zs{n}}{n^{\epsilon}}
=0.
\end{equation*}
}

\vspace{2mm}

\noindent $\C_\zs{2})$ {\it Assume that the set $\Lambda$ is such that for any $\epsilon>0$
\begin{equation}\label{sec:Mo.8+1++2}
\lim_\zs{n\to\infty}\frac{\nu_n}{n^{\epsilon}}=0
\quad\mbox{and}\quad
\lim_\zs{n\to\infty}\,\frac{\vert\Lambda\vert_\zs{n}}{n^{1/2+\epsilon}}
=0\,.
\end{equation}
}

\vspace{2mm}

We note that Theorem \ref{sec:Mo.Th.1} and Proposition \ref{sec:Mo.Prop.1} directly
imply the following inequality.
\begin{theorem}\label{sec:Mo.Th.2}
If the conditions  $\C_\zs{1})$ -- $\C_\zs{2})$ hold for the
distribution $Q$ of the process $\xi$ in \eqref{sec:In.1}, then
for any $n\geq 2$ and $0<\delta<1/2$, the robust risk \eqref{sec:risks_11_0} of estimate \eqref{sec:Mo.6} for
continuously differentiable function
$S$
satisfies the oracle inequality
\begin{equation}\label{sec:Mo.10+1-2}
\cR^{*}_\zs{n}
(S^*,S)\,\le\, \frac{1+5\delta}{1-\delta}
\min_\zs{\lambda\in\Lambda} \cR^{*}_\zs{n}\,(S^*_\zs{\lambda},S)
+\frac{B_\zs{n}(1+\|\dot{S}\|^2)}{n \delta}\,,
\end{equation}
where the term $B_\zs{n}$ is
independent of $S$ and for any $\epsilon>0$
\begin{equation*}\label{sec:Mo.11-2}
\lim_\zs{n\to\infty}\,\frac{B_\zs{n}}{n^{\epsilon}}
=0\,.
\end{equation*}
\end{theorem}

\begin{remark}\label{Re.sec:ModSel_010}
Note that  sharp oracle inequalities
similar to
\eqref{sec:Mo.10} and \eqref{sec:Mo.10+1-2}
was obtained earlier by
Konev and Pergamenshchikov
 in
 \cite{KonevPergamenshchikov2009a, KonevPergamenshchikov2012, KonevPergamenshchikov2015}
 for model selection procedures based on the weighted  least squares estimates \eqref{sec:Imp.5}.
Unfortunately, we can not use such oracle inequalities for the model selection procedures,
based on the weighted shrinkage estimates \eqref{sec:Imp.11}
since
they  depend non linearly on the coefficients $\lambda$. This is main technical difficulty which
doesn't allow us to use the obtained oracle inequalities. Moreover, in all these papers
the oracle inequalities are obtained under condition that
 the L\'evy measure is finite. The inequalities
 \eqref{sec:Mo.10} and \eqref{sec:Mo.10+1-2}
 are obtained without conditions on the impulse noises.
\end{remark}

Now we specify the weight coefficients $(\lambda(j))_\zs{j\ge
1}$ in the way proposed in \cite{GaltchoukPergamenshchikov2009a}
 for a heteroscedastic regression
model in discrete time. Consider a numerical grid of the form
\begin{equation*}\label{sec:Imp.7}
\cA_\zs{n}=\{1,\ldots,k_\zs{n}\}\times\{r_1,\ldots,r_m\}\,,
\end{equation*}
where  $r_i=i\rho_n$ and $m=[1/\rho_n^2]$. Both
parameters $k_n\ge 1$ and $0<\rho_n\le 1$ are the
functions of $n$ such that for any $\epsilon>0$
\begin{equation}\label{sec:Imp.8}
\left\{
\begin{array}{ll}
&\lim_\zs{n\to\infty}\,k_n=+\infty\,,
\quad
\lim_\zs{n\to\infty}\,\dfrac{k_n}{\ln n}=0\,,\\[6mm]
&\lim_\zs{n\to\infty}\rho_n=0
\quad\mbox{and}\quad
\lim_\zs{n\to\infty}\,n^{\epsilon}\rho_n\,=+\infty\, .
\end{array}
\right.
\end{equation}
One can take, for
example,
$$
\rho_n=\frac{1}{\ln (n+1)}
\quad\mbox{and}\quad
k_n=\sqrt{\ln (n+1)}\,.
$$
For each $\alpha=(\beta,r)\in\cA_\zs{n}$ we introduce the weight
sequence $\lambda_\zs{\alpha}=(\lambda_\zs{j}(\alpha))_\zs{j\ge 1}$
as
\begin{equation}\label{sec:Imp.9}
\lambda_\zs{j}(\alpha)=\Chi_\zs{\{1\le j\le d\}}+
\left(1-(j/\omega_\alpha)^\beta\right)\, \Chi_\zs{\{ d<j\le
\omega_\alpha\}}
\end{equation}
where $d=d(\alpha)=\left[\omega_\zs{\alpha}/\ln (n+1)\right]$,
$$
\omega_\zs{\alpha}=\left(\tau_\zs{\beta}\,r\,v_\zs{n}\right)^{1/(2\beta+1)}\,,
\quad
\tau_\zs{\beta}=\frac{(\beta+1)(2\beta+1)}{\pi^{2\beta}\beta} \quad\mbox{and}\quad v_\zs{n}=
\frac{n}{\overline{\sigma}_n}
\,.
$$
We set
\begin{equation}\label{sec:Imp.10_Lambda}
\Lambda\,=\,\{\lambda(\alpha)\,,\,\alpha\in\cA_\zs{n}\}\,.
\end{equation}
It will be noted that in this case $\nu_n=k_n m$.  Therefore, the conditions
\eqref{sec:Imp.8} imply the first limit equality in \eqref{sec:Mo.8+1++2}. Moreover,
in view of the definition \eqref{sec:Imp.9}
and taking into account that $\tau_\zs{\beta}\le 1$ for $\beta\ge 1$
 the function $L(\lambda)$ defined in
\eqref{sec:Imp.6+0}
can be estimated
for any $\lambda\in\Lambda$
as
$$
\max_\zs{\lambda\in\Lambda}\,
L(\lambda)\le \max_\zs{\lambda\in\Lambda}\,\omega_\zs{\alpha}\le v^{1/3}_\zs{n} \rho^{-1/3}_n
\,.
$$
Therefore, using here
the conditions
\eqref{sec:Ex.01-2}
and
\eqref{sec:Imp.8}
we get the last  limit  in \eqref{sec:Mo.8+1++2}, i.e.
the condition $\C_2)$ holds for the set
 $\Lambda$ defined in \eqref{sec:Imp.10_Lambda}.

\begin{remark}\label{Re.sec:ModSel+++000}
It will be observed that the specific form of weights
\eqref{sec:Imp.9} was proposed by Pinsker \cite{Pinsker1981} for the
filtration problem with known smoothness of the regression
function observed with an additive gaussian white noise in
continuous time. Nussbaum
\cite{Nussbaum1985}
 used such weights for the
gaussian regression estimation problem in discrete time.
\end{remark}

\bigskip
\section{Monte Carlo simulations}\label{sec:Sim}

In this section we give the results of numerical simulations to assess the performance and improvement of
the proposed model selection procedure \eqref{sec:Mo.6}.
We simulate the model \eqref{sec:In.1} with
$1$-periodic function $S$ of the form
\begin{equation}\label{sec:Sim_Sign_11}
S(t)=t\,\sin(2\pi t)+t^2(1-t)\cos(4\pi t)
\end{equation}
on $[0,\,1]$ and the L\'evy noise process $\xi_\zs{t}$ is defined as
$$
\xi_\zs{t}=0.5\, w_\zs{t}+0.5\, z_\zs{t}\,.
$$
Here $z_\zs{t}$ is a compound Poisson process with intensity $\lambda=\Pi(x^2)=1$ and a Gaussian
$\cN(0,\,1)$ sequence $(Y_\zs{j})_\zs{j\ge1}$ (see, for example, \cite{KonevPergamenshchikov2015}).

We use the model selection procedure   \eqref{sec:Mo.6} with the weights \eqref{sec:Imp.9} in which
$k_n=100+\sqrt{\ln (n+1)}$, $r_i=i/\ln (n+1)$, $m=[\ln^2 (n+1)]$, $\overline{\sigma}_n=0.5$ and $\delta=(3+\ln n)^{-2}$.
We define the empirical risk as
$$
\cR(S^*,\,S)=\frac{1}{p}\sum_\zs{j=1}^p \wh{\E}\left(S_n^*(t_j)-S(t_j)\right)^2\,,
$$
$$
 \wh{\E}\left(S_n^*(\cdot)-S(\cdot)\right)^2=
\frac{1}{N}\sum_\zs{l=1}^N \left(S_{n,l}^*(\cdot)-S(\cdot)\right)^2\,,
$$
where the observation frequency $p=100001$ and the expectations was taken as an
average over $N = 1000$ replications.


\begin{table}[h]
\label{Tab1}
\caption{The sample quadratic risks for different optimal $\lambda$}
\begin{center}
\begin{tabular}{|c|c|c|c|c|}
  \hline
   $n$              & 100 & 200 & 500 & 1000 \\ \hline
  $\cR(S^*_\zs{\lambda^*},\,S)$    & 0.0118 & 0.0089 & 0.0031 & 0.0009 \\ \hline
  $\cR(\wh{S}_\zs{\wh{\lambda}},\,S)$ & 0.0509 & 0.0203 & 0.0103 & 0.0064 \\ \hline
  $\cR(\wh{S}_\zs{\wh{\lambda}},\,S)/\cR(S^*_\zs{\lambda^*},\,S)$ & 4.3 & 2.3 & 3.3 & 7.1 \\
  \hline
\end{tabular}
\end{center}
\end{table}

\begin{table}[h]
\label{Tab2}
\caption{The sample quadratic risks for the same optimal $\wh{\lambda}$}
\begin{center}
\begin{tabular}{|c|c|c|c|c|}
  \hline
   $n$              & 100 & 200 & 500 & 1000 \\ \hline
  $\cR(S^*_\zs{\wh{\lambda}},\,S)$    & 0.0237 & 0.0103 & 0.0041 & 0.0011 \\ \hline
  $\cR(\wh{S}_\zs{\wh{\lambda}},\,S)$ & 0.0509 & 0.0203 & 0.0103 & 0.0064 \\ \hline
  $\cR(\wh{S}_\zs{\wh{\lambda}},\,S)/\cR(S^*_\zs{\wh{\lambda}},\,S)$ & 2.1 & 2.2 & 2.5 & 5.8 \\
  \hline
\end{tabular}
\end{center}
\end{table}

Table 1 gives the values for the sample risks of the improved estimate \eqref{sec:Mo.6}
and the model selection procedure based on the weighted LSE (3.15) from \cite{KonevPergamenshchikov2012} for different numbers
of observation period $n$. Table 2 gives the values for the sample risks of the the model selection procedure based on the weighted LSE (3.15) from \cite{KonevPergamenshchikov2012} and it's improved version for different numbers
of observation period $n$.


\begin{figure}[h!]
\label{fig1}
\centering
    \includegraphics[width=0.7\textwidth]{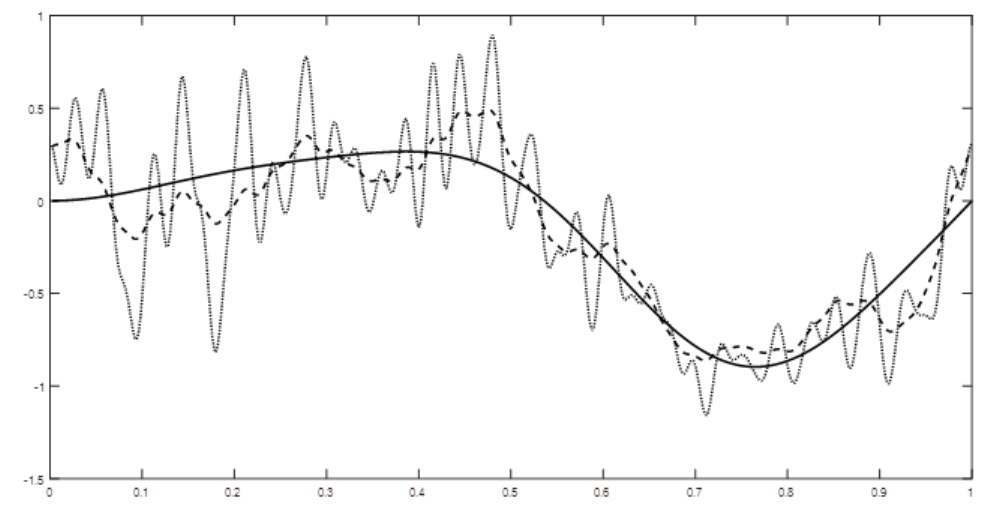}
  \caption{Behavior of the regression function and its estimates for $n=100$.}
\end{figure}

\newpage

\begin{figure}[h!]
\label{fig2}
\centering
    \includegraphics[width=0.7\textwidth]{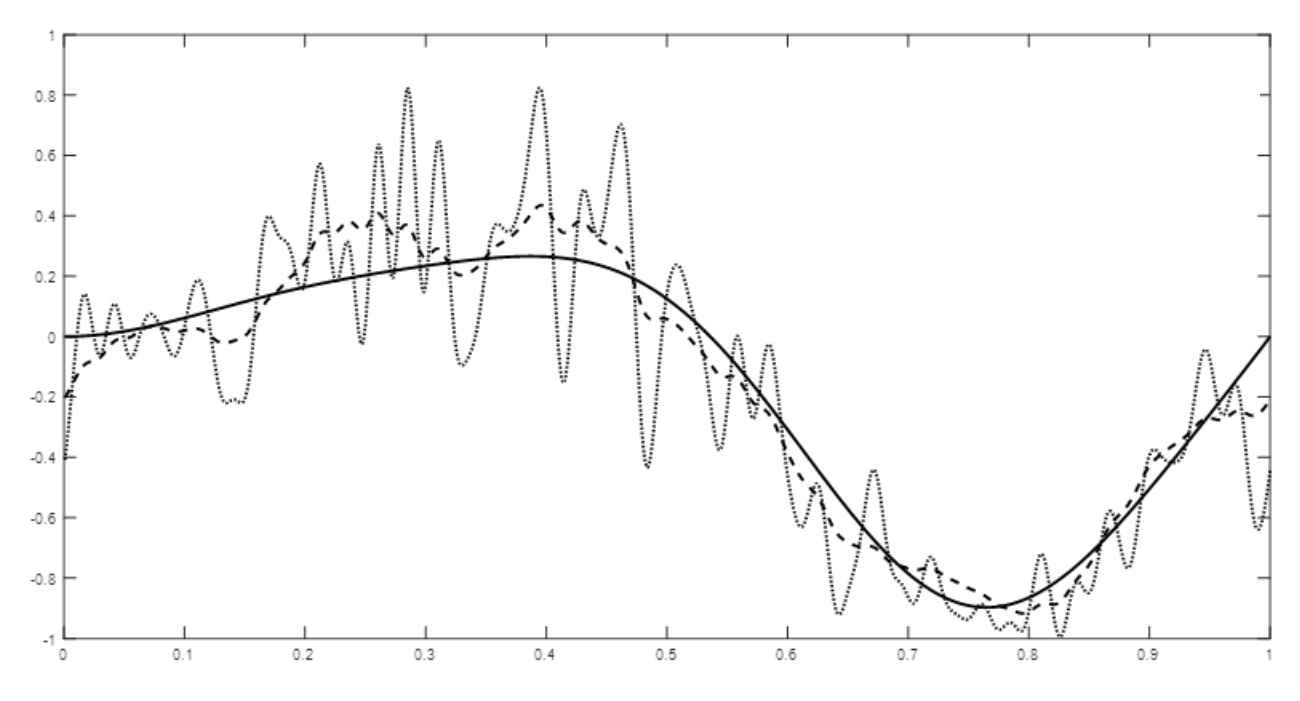}
  \caption{Behavior of the regression function and its estimates for $n=200$.}
\end{figure}


\begin{figure}[h!]
\label{fig3}
\centering
    \includegraphics[width=0.7\textwidth]{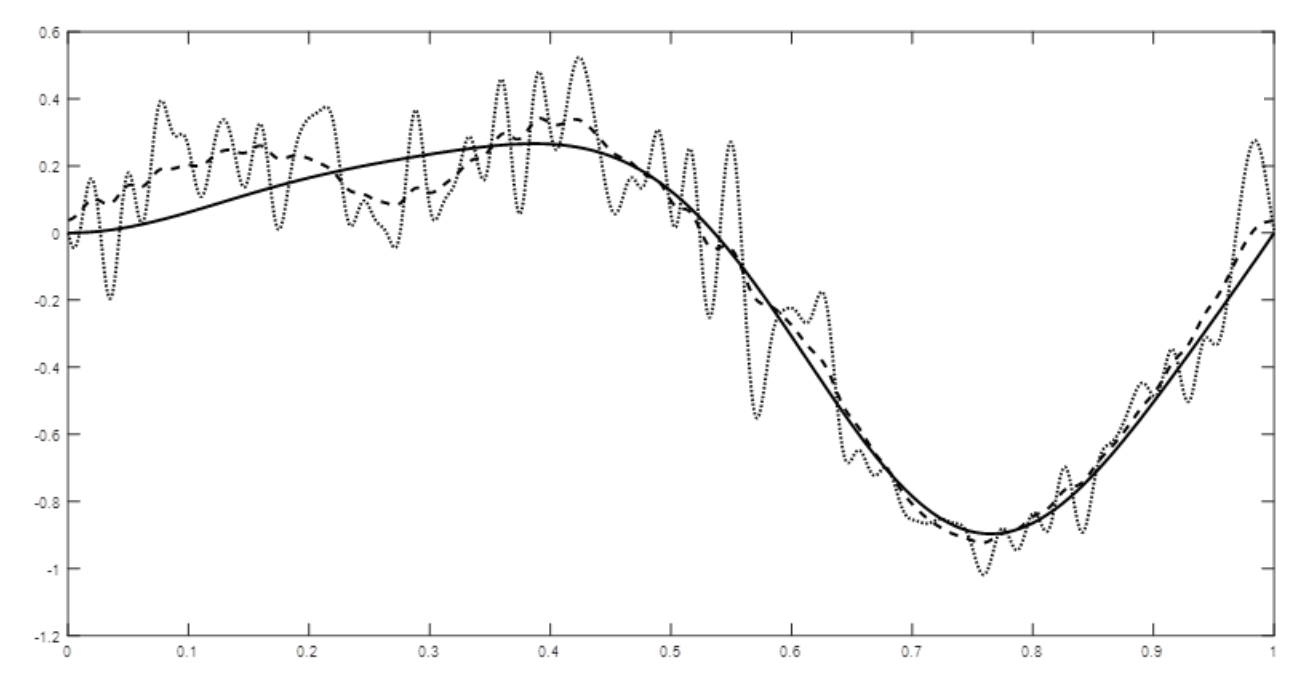}
  \caption{Behavior of the regression function and its estimates for $n=500$.}
\end{figure}


\begin{figure}[h!]\label{fig4}
\centering
 \includegraphics[width=0.7\textwidth]{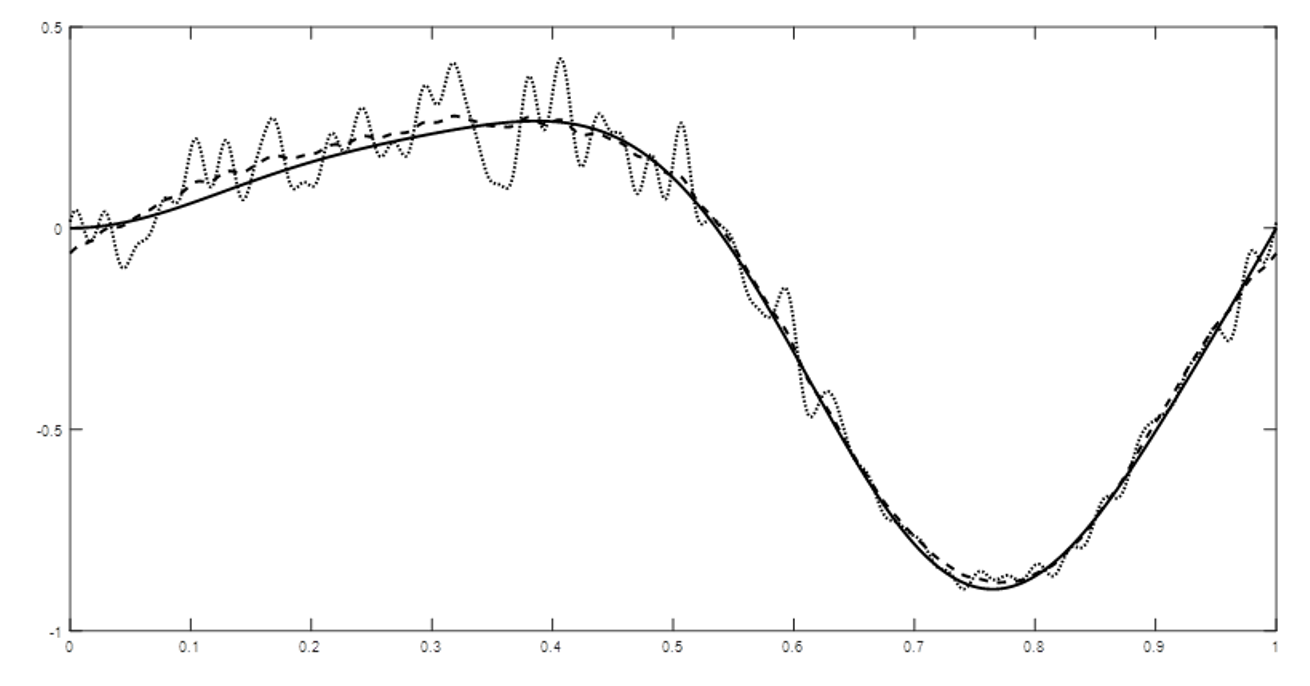}
  \caption{Behavior of the regression function and its estimates for $n=1000$.}
\end{figure}

\bigskip

\begin{remark}
Figures 1--4 show the behavior of
the procedures  \eqref{sec:Imp.5} and \eqref{sec:Mo.6}   depending
on the values of observation periods $n$.
The bold line is the function \eqref{sec:Sim_Sign_11},
the continuous  line is the  model selection
procedure based on the
least squares estimators $\wh{S}$
 and the dashed line is the improved model selection procedure $S^*$.
 From the Table 2 for the same $\lambda$ with various observations
numbers $n$ we can conclude that
theoretical result on the improvement effect \eqref{sec:Imp.11+1} is confirmed
by the numerical simulations.
Moreover, for the proposed shrinkage procedure, Table 1 and Figures 1--4,
we can conclude that the benefit is considerable for non large  $n$.
\end{remark}

\bigskip

\section{Asymptotic efficiency}\label{sec:Ae}

In order to study the asymptotic efficiency we define the following functional Sobolev ball
\begin{equation*}\label{sec:Ae.1}
W_\zs{k,r}=\{f\in\C^{k}_\zs{p}[0,1]\,:\,
\sum_\zs{j=0}^k\,\|f^{(j)}\|^2\le r\}\,,
 \end{equation*}
where $r>0$ and $k\ge 1$ are
some unknown parameters, $\C^{k}_\zs{p}[0,1]$ is the space of
 $k$ times differentiable $1$ - periodic $\bbr\to\bbr$ functions
 such that for any $0\le i \le k-1$
$$
f^{(i)}(0)=f^{(i)}(1)
\,.
$$
 In order to formulate our asymptotic results
we define the Pinsker constant which gives the lower bound for normalized asymptotic risks
\begin{equation}\label{sec:Ae.3}
l_\zs{k}(r)\,=\,((1+2k) r)^{1/(2k+1)}\,
\left(\frac{k}{\pi (k+1)}\right)^{2k/(2k+1)}
\,.
\end{equation}

It is well known that for any $S\in W_\zs{k,r}$
 the optimal rate of convergence is
$n^{-2k/(2k+1)}$ (see, for example, \cite{Pinsker1981, Nussbaum1985}).
On the basis of the model selection procedure
we construct the adaptive
procedure $S^*$ for which we obtain the following asymptotic upper bound
for the quadratic risk, i.e.
 we show that the parameter \eqref{sec:Ae.3} gives a lower bound
for the asymptotic normalized risks.
To this end we denote by $\Sigma_\zs{n}$ the set of all estimators $\wh{S}_\zs{n}$ of $S$ measurable with respect to
the process
\eqref{sec:In.1}, i.e. measurable with respect to $\sigma$-field
$\sigma\{y_\zs{t}\,,\,0\le t\le n\}$.

\begin{theorem}\label{Th.sec: Ae.2}
The robust risk \eqref{sec:risks_11_0}
admits the following asymptotic lower bound
 \begin{equation*}\label{sec:Ae.05}
\liminf_\zs{n\to\infty}\,
\inf_\zs{\wh{S}_\zs{n}\in\Sigma_\zs{n}}
\,v_\zs{n}^{2k/(2k+1)}
\sup_\zs{S\in W_\zs{k,r}}\,\cR^{*}_\zs{n}(\wh{S}_\zs{n},S)
\,
\ge l_\zs{k}(r)  \,.
 \end{equation*}
\end{theorem}

\noindent
We show that this lower bound is sharp in the following sense.

\begin{theorem}\label{Th.sec: Ae.1}
The quadratic risk \eqref{sec:risks_11_0} for the
 estimating procedure $S^{*}$ has the following asymptotic upper bound
 \begin{equation*}\label{sec:Ae.4}
\limsup_\zs{n\to\infty}\,v_\zs{n}^{2k/(2k+1)}
\sup_\zs{S\in W_\zs{k,r}}\,\cR^{*}_\zs{n}(S^{*},S)
\,
\le l_\zs{k}(r)
\,.
 \end{equation*}
\end{theorem}

\noindent
It is clear that Theorem \ref{Th.sec: Ae.1}
 and Theorem \ref{Th.sec: Ae.2}
imply
\begin{corollary}\label{Co.sec: Ae.1}
The model selection procedure $S^{*}$
is asymptotically efficient, i.e.
\begin{equation}\label{sec:Ae.5}
\lim_{n\to\infty}\,(v_\zs{n})^{\frac{2k}{2k+1}}\,
\sup_\zs{S\in W_\zs{k,r}}\,\cR^{*}_\zs{n}(S^{*},S)\,
= l_\zs{k}(r)
\,.
\end{equation}
\end{corollary}

\begin{remark}\label{Re.2.2}
Note that the equality
\eqref{sec:Ae.5} implies that the parameter
\eqref{sec:Ae.3} is the Pinsker constant
in this case (cf. \cite{Pinsker1981}).
\end{remark}

\begin{remark}\label{Re.2.3}
It should be noted that the equality \eqref{sec:Ae.5}
means that the robust efficiency holds with
the  convergence rate
$$
(v_\zs{n})^{\frac{2k}{2k+1}}
\,.
$$
It is well known that for the simple risks
 the optimal (minimax) estimation convergence rate
for the functions from the set  $W_\zs{k,r}$
 is $n^{2k/(2k+1)}$ (see, for example, \cite{Pinsker1981, Nussbaum1985, IbragimovKhasminskii1981}).
So,
 if  the distribution upper bound  $\overline{\sigma}_n\to 0$  as $n\to\infty$
 we obtain the more rapid rate, and
 if $\overline{\sigma}_n\to \infty$  as $n\to\infty$
we obtain the more slow rate. In the case when $\overline{\sigma}_n$ is constant the robust rate is the same as the classical non robust convergence rate.
\end{remark}

\begin{remark}
The property
\eqref{sec:Ae.5} means that the model selection procedure
\eqref{sec:Mo.6} asymptotically has the same efficiency property as
the LSE model selection (see, \cite{GaltchoukPergamenshchikov2009b, KonevPergamenshchikov2009b}).
So, it means that the proposed shrinkage method non-asymptotically has benefit
with respect to LSE and asymptotically the shrinkage methods keep the efficiency property.
\end{remark}

\bigskip

\section{Stochastic calculus  for the L\'evy processes}\label{sec:Stc}

In this section we study the process \eqref{sec:In.1+1}.
First we recall  the Novikov  inequalities, \cite{Novikov1975}, also referred to as the Bichteler--Jacod inequalities, see
\cite{BichtelerJacod1983, MarinelliRockner2014},  providing bounds of the moments of the supremum of purely discontinuous local martingales for $p\ge 2$
and for any $n\ge 1$
\begin{equation}
\label{Novikov++}
\E\sup_\zs{0\le t\le n}|\Upsilon*(\mu-\wt{\mu})_\zs{t}|^{p}\le C_\zs{p}
 \left(
 \E\,\big (|\Upsilon|^{2}*\wt{\mu}_\zs{n}\big)^{p/2}
 +
 \E\,\big (|\Upsilon|^{p}*\wt{\mu}_\zs{n}\big)
\right)\,,
\end{equation}
where $C_\zs{p}$ is some positive constant.  Further for any two functions $f$ and $g$ from $\L_\zs{2}[0,t]$  with $t>0$ we
use the following notations
\begin{equation}
\label{product_norm}
(f,g)_\zs{t}=
\int^{t}_\zs{0}
f(s)g(s)
\d s
\quad\mbox{and}\quad
\Vert f\Vert^{2}_\zs{t}=
\int^{t}_\zs{0}
f^{2}(s)
\d s
\,.
\end{equation}

\begin{proposition}\label{Pr.sec:Stc.1}
For any nonrandom function $f$ and $g$ from $\L_\zs{2}[0,t]$
\begin{equation}\label{sec:Stc.1}
\E\, I_\zs{t}(f)I_\zs{t}(g)=
\sigma\,
(f,g)_\zs{t}\,,
\end{equation}
where the noise variance $\sigma$ is given in \eqref{sec:Ex.01-1}.
\end{proposition}

\noindent
 Now we set
\begin{equation}\label{sec:Stc.6-00_wt_I}
\wt{I}_\zs{t}(f)=I^{2}_\zs{t}(f)
-\E\,I^{2}_\zs{t}(f)
\quad\mbox{and}\quad
\wt{M}_\zs{t}(f)=\M^{f,f}_\zs{t}
\,.
\end{equation}

\noindent
For any  $[0,n]\to\bbr$  function $f$ we introduce the following uniform norm
$$
\|f\|_\zs{n}=\sup_\zs{0\le t\le n}|f(t)|\,.
$$

\begin{proposition}
\label{Pr.sec:Stc.3}
Let $f$ and $g$ be two
borel  $[0,n]\to\bbr$ functions such that
$\|f\|_\zs{n}\le \overline{\phi}_\zs{n}$ and $\|g\|_\zs{n}\le \overline{\phi}_\zs{n}$.
Then for any $0<t\le n$
\begin{equation}\label{sec:Stc.9+1_UpBnd}
\left\vert \E\,\wt{I}_\zs{t}(f)\, \wt{I}_\zs{t}(g)\,
\right\vert\le
\sigma^{2}\,
\left(
2(f,g)^{2}_\zs{t}
+
\overline{\phi}^{4}_\zs{n}\,\Pi(x^4)t
\right)\,.
\end{equation}
\end{proposition}

\proof Using \eqref{sec:Stc.2} with $f=g$ we can obtain that
the process $(\wt{I}_\zs{t})_\zs{t\ge 0}$ satisfies the following
stochastic equation
\begin{equation*}\label{sec:Stc.6-01}
\d \wt{I}_\zs{t}(f)
=
\d\wt{M}_\zs{t}(f)\,,
\quad
\wt{I}_\zs{0}(f)=0\,.
\end{equation*}
Note that from the definition of $M_\zs{t}(f,f)$ in
\eqref{sec:Stc.2}
we can represent $\wt{I}_\zs{t}(f)$ as
\begin{equation}\label{sec:Stc.6-01_reprnt}
\wt{I}_\zs{t}(f)=\wt{I}^{c}_\zs{t}(f)+\wt{I}^{d}_\zs{t}(f)
\,,
\end{equation}
where $\wt{I}^{c}_\zs{t}(f)=2\sigma_\zs{1}\int^{t}_\zs{0}\,I_\zs{s}(f)f(s)\d w_\zs{s}$
and
$$
\wt{I}^{d}_\zs{t}(f)=2\sigma_\zs{2}\int^{t}_\zs{0}\,I_\zs{s-}(f)f(s)\d z_\zs{s}
+
\sigma^{2}_\zs{2}\,
\int^{t}_\zs{0}
f^{2}(s)\,
\d m_\zs{s}
\,.
$$
Moreover, by the Ito formula
\begin{align*}
\wt{I}_\zs{t}(f)\,\wt{I}_\zs{t}(g)&=
\int^{t}_\zs{0}
\wt{I}_\zs{s-}(f)\,\d\wt{I}_\zs{s}(g)
+
\int^{t}_\zs{0}
\wt{I}_\zs{s-}(g)\,\d\wt{I}_\zs{s}(f)\\[2mm]
&+4\sigma^{2}_\zs{1}
\int^{t}_\zs{0}
f(s)\,g(s)I_\zs{s}(f)I_\zs{s}(g)\d s
+
\check{J}^{f,g}_\zs{t}\,,
\end{align*}
where $\check{J}^{f,g}_\zs{t}=\sum_\zs{0<s\le t}\,\Delta \wt{I}^{d}_\zs{s}(f)\,\Delta\wt{I}^{d}_\zs{s}(g)$.
Using
the last in condition
\eqref{sec:Ex.1-00_mPi}
and
the inequality
\eqref{Novikov++}
we can obtain that
for any bounded measurable $[0,n]\to\bbr$ function $h$
\begin{equation}\label{sec:Upper_bound_h}
\sup_\zs{0\le t\le n}\,
\E
\,
\left(
\int^{t}_\zs{0}\,h(s)\,\d z_\zs{s}
\right)^{6}
\,
<\infty\,.
\end{equation}
From this and the H\"older inequality we obtain that
$$
\sup_\zs{0\le t\le n}\,
\E\,
I^{4}_\zs{t}(g)
\,
I^{2}_\zs{t}(f)
\,
\le
\,
\sup_\zs{0\le t\le n}\,
\left(\E\,
I^{6}_\zs{t}(g)
\right)^{2/3}
\,
\left(\E\,
I^{6}_\zs{t}(f)
\right)^{1/3}
\,<\infty\,.
$$
Therefore, in view of Proposition \eqref{Pr.sec:Stc.1}
$$
\E\,\wt{I}_\zs{t}(f)\,\wt{I}_\zs{t}(g)=
4\sigma\sigma^{2}_\zs{1}\int^{t}_\zs{0}\,f(s)\,g(s)(f,g)_\zs{s}\d s
+\E\,\check{J}^{f,g}_\zs{t}
=
2\sigma\sigma^{2}_\zs{1}
(f,g)^{2}_\zs{t}
+\E\,\check{J}^{f,g}_\zs{t}
\,.
$$
From the definition of the discrete part of $\wt{I}^{d}_\zs{t}(f)$ in \eqref{sec:Stc.6-01_reprnt}
we can represent the jumps term $\check{J}^{f,g}_\zs{t}$ as
\begin{equation}
\label{sec:jumps_reprent_t}
\check{J}^{f,g}_\zs{t}=4\sigma^{2}_\zs{2} \check{J}^{f,g}_\zs{1,t}
+2\sigma^{3}_\zs{2} \check{J}^{f,g}_\zs{2,t}
+
\sigma^{4}_\zs{2} \check{J}^{f,g}_\zs{3,t}
\,,
\end{equation}
where $\check{J}^{f,g}_\zs{1,t}=\sum_\zs{0<s\le t}
I_\zs{s-}(f)I_\zs{s-}(g) f(s)g(s)\,\left( \Delta z_\zs{s}\right)^{2}$,
$$
\check{J}^{f,g}_\zs{2,t}=\sum_\zs{0<s\le t}
\left(
I_\zs{s-}(f) f(s) g^{2}(s)\,
+
I_\zs{s-}(g)( f^{2}(s)\,
\right)
\left( \Delta z_\zs{s}\right)^{3}
$$
and
$\check{J}^{f,g}_\zs{3,t}=\sum_\zs{0<s\le t}\, f^{2}(s)g^{2}(s)\,\left( \Delta z_\zs{s}\right)^{4}$.
In view of Proposition \ref{Pr.sec:Stc.1}
and the upper bound \eqref{sec:Upper_bound_h} and taking into account that $\Pi(x^{2})=1$
we calculate
$$
\E\,\check{J}^{f,g}_\zs{1,t}=
\int^{t}_\zs{0}\,\E\,I_\zs{s}(f)I_\zs{s}(g) f(s)g(s)\d s
=\sigma\int^{t}_\zs{0}\,(f,g)_\zs{s} f(s)g(s)\d s
=
\frac{\sigma}{2}\,
(f,g)^{2}_\zs{t}
\,.
$$
Similarly, we obtain that
$$
\E\,\check{J}^{f,g}_\zs{2,t}=\Pi(x^{3})
\int^{t}_\zs{0}\,f(s)g(s)
\left(g(s)
\E\,I_\zs{s}(f)
+
f(s)
\E\,I_\zs{s}(g)
\right)
\d s =0\,.
$$
and
$\E\check{J}^{f,g}_\zs{3,t}=\Pi(x^{4})\int^{t}_\zs{0}\, f^{2}(s)g^{2}(s)\,\d s$. So,
$$
\E\,\check{J}^{f,g}_\zs{t}
=
2 \sigma^{2}_\zs{2} \sigma (f,g)^{2}_\zs{t}
+
\sigma^{4}_\zs{2} \Pi(x^{4})\,
(f^{2},g^{2})_\zs{t}
\,,
$$
and, therefore,
$$
\E\,\wt{I}_\zs{t}(f)\,\wt{I}_\zs{t}(g)=
2\sigma^{2}\,
(f,g)^{2}_\zs{t}
+
\sigma^{4}_\zs{2}
 \Pi(x^{4})\,
(f^{2},g^{2})_\zs{t}
\,.
$$
Taking into account here that $\sigma^{4}_\zs{2} \le \sigma^{2}$
 and the conditions of the proposition we obtain the upper bound \eqref{sec:Stc.9+1_UpBnd}.
Hence Proposition \ref{Pr.sec:Stc.3}. \endproof

\noindent Now for any $y\in\bbr^{n}$ we define the following function
\begin{equation}\label{sec:Stc.9-00+1_Ix}
\overline{I}_\zs{t}(y)=
\sum^{n}_\zs{j=1}\,y_\zs{j}\,\wt{I}_\zs{t}(\phi_\zs{j})\,.
\end{equation}

\noindent
For this we show the following property.
\begin{proposition}
\label{Pr.sec:Stc.4}
For any $n\ge 1$
\begin{equation}\label{sec:Stc.10+1}
\sup_\zs{y\in\bbr^{n}\,,\,\Vert y\Vert\le 1}
\E\,\overline{I}^{2}_\zs{n}(y)
\le\,(2+\overline{\phi}_\zs{n}^4\Pi(x^4))
\sigma^2\,
n^{2}
\,.
\end{equation}
\end{proposition}

\proof From  Proposition \ref{Pr.sec:Stc.3} it follows that
\begin{align*}
\E\,\overline{I}^{2}_\zs{n}(y)&=
\sum^{n}_\zs{i,j=1}y_\zs{i}\,y_\zs{j}\,\E\,
\wt{I}_\zs{n}(\phi_\zs{i}) \wt{I}_\zs{n}(\phi_\zs{j})\\[2mm]
&\le 2\sigma^2
\sum^{n}_\zs{i,j=1}\vert y_\zs{i}\vert\,\vert y_\zs{j}\vert
\left(\phi_\zs{i}\,,\,\phi_\zs{j}\right)^{2}_\zs{n}
+n\overline{\phi}^{4}_\zs{n}
\sigma^2
\Pi(x^4)
\left(\sum^{n}_\zs{i=1}\vert y_\zs{i}\vert\right)^{2}
\,.
\end{align*}
Taking into account here that the functions $(\phi_\zs{j})_\zs{j\ge 1}$ are orthonormal, and the fact that
$\left(\sum^{n}_\zs{i=1}\vert y_\zs{i}\vert\right)^{2}\le n$, we obtain the bound \eqref{sec:Stc.10+1}.
\endproof

\section{The van Trees inequality for the L\'evy processes}\label{sec:VanTrees}

In this section we consider the following continuous time
 parametric regression model \eqref{sec:In.1}
with the  function $S$ defined as
\begin{equation*}\label{sec:VanTrees_1}
S(t,\theta)=
\sum^{d}_\zs{i=1}\,\theta_\zs{i}\,\psi_\zs{i}(t)\,,
\end{equation*}
with the unknown parameters
$\theta=(\theta_\zs{1},\ldots,\theta_\zs{d})'$. Here we assume that the functions $(\psi_i)_\zs{1\le i\le d}$ are $1$-periodic and orthogonal functions.

Let us denote by $\nu_\zs{\xi}$ the distribution of the process $(\xi_\zs{t})_\zs{0\le t\le n}$
on the Skorokhod space $\D[0,n]$.  From Proposition \ref{Pr.sec:App.1++} it follows that in this space for any parameters $\theta\in\bbr^d$, the distribution
$\P_\zs{\theta}$
of the process \eqref{sec:In.1}
 is absolutely continuous with respect to the
$\nu_\zs{\xi}$
 and the corresponding Radon-Nikodym derivative,
for any function $x=(x_\zs{t})_\zs{0\le t\le n}$ from $\D[0,n]$,
 is defined as
\begin{equation}\label{sec:App.7}
f(x,\theta)=
\frac{\d\P_\zs{\theta}}{\d\nu_\zs{\xi}}(x)=
\exp\left\{\int^{n}_\zs{0}\,\frac{S(t,\theta)}{\sigma^{2}_\zs{1}}\,
\d x^{c}_\zs{t}
-\,\int^{n}_\zs{0}\,
\frac{S^{2}(t,\theta)}{2\sigma^{2}_\zs{1}}\,
\d t
\right\}
\,,
\end{equation}
where
$$
x^{c}_\zs{t}=
x_\zs{t}
-
\int^{t}_\zs{0}\,\int_\zs{\bbr}\,v\,\left(
\mu_\zs{x}(\d s\,,\d v)
-
\Pi(\d v)\d s
\right)
$$
and for any measurable set $A$ in $\bbr$ with $0\notin A$
$$
\mu_\zs{x}([0,t]\times A)=\sum_\zs{0\le s\le t}\,
\Chi_\zs{\{\Delta\xi_\zs{s}\in \sigma_\zs{2} A\}}
\,.
$$

\noindent
Let $U$ be a prior density on $\bbr^d$ having
the following form:
$$
U(\theta)=U(\theta_1,\ldots,\theta_d)=\prod_{j=1}^d u_\zs{j}(\theta_\zs{j})\,,
$$
where $u_\zs{j}$ is some continuously differentiable density in $\bbr$.
Moreover, let $g(\theta)$ be a continuously differentiable $\bbr^d\to \bbr$ function such that,
for each $1\le j\le d$,
\begin{equation}\label{sec:App.8}
\lim_\zs{|\theta_\zs{j}|\to\infty}\,
g(\theta)\,u_\zs{j}(\theta_\zs{j})=0
\quad\mbox{and}\quad
\int_\zs{\bbr^d}\,|g^{\prime}_\zs{j}(\theta)|\,U(\theta)\,\d \theta
<\infty\,,
\end{equation}
where $g^{\prime}_\zs{j}(\theta)=\,\partial g(\theta)/\partial\theta_\zs{j}$.
For any $\cB(\cX)\times\cB(\bbr^d)$-measurable integrable function $H=H(x,\theta)$ we denote
\begin{align*}
\wt{\E}\,H&=\int_{\bbr^d}\,
\int_\zs{\cX}\,H(x,\theta)\,\d \P_\zs{\theta}\,U(\theta) \d \theta\\[2mm]
&=
\int_{\bbr^d}\,\int_\zs{\cX}\,
H(x,\theta)\,f(x,\theta)\,U(\theta)\d \nu_\zs{\xi}(x)\, \d \theta\,,
\end{align*}
where $\cX=\D[0,n]$.

\begin{lemma}\label{Le.sec:App.3}
For any $\cF^y_n$-measurable square integrable function $\wh{g}_\zs{n}$
 and for any $1\le j\le d$, the following inequality holds
$$
\wt{\E}(\wh{g}_\zs{n}-g(\theta))^2\ge
\frac{\eta^2_\zs{j}}{n\Vert \psi_\zs{j}\Vert^{2}\sigma^{-2}_\zs{1}+I_\zs{j}}\,,
$$
where
$$
\eta_\zs{j}=\int_\zs{\bbr^d}\,g^{\prime}_\zs{j}(\theta)\,U(\theta)\,\d \theta
\quad\mbox{and}\quad
I_\zs{j}=\int_\zs{\bbr}\,\frac{\dot{u}^2_\zs{j}(z)}{u_\zs{j}(z)}\,\d z\,.
$$
\end{lemma}
\noindent {\bf  Proof.}
First of all note that, the density \eqref{sec:App.7}
on the process $\xi$
is bounded
with respect to $\theta_\zs{j}\in\bbr$ and for any $1\le j\le d$
$$
\limsup_\zs{|\theta_\zs{j}|\to\infty}\,f(\xi,\theta)\,=\,0\,
\quad\quad\mbox{a.s.}
$$
Now, we set
$$
\wt{U}_\zs{j}=\wt{U}_\zs{j}(x,\theta)=
\frac{\partial\,(f(x,\theta)U(\theta))/\partial\theta_\zs{j}}{f(x,\theta)U(\theta)}
 \,.
$$
Taking into account the condition \eqref{sec:App.8} and
integrating by parts yield
\begin{align*}
\wt{\E}&\left((\wh{g}_\zs{n}-g(\theta))\wt{U}_\zs{j}\right)
=\int_{\cX\times\bbr^d}\,(\wh{g}_\zs{n}(x)-g(\theta))\frac{\partial}{\partial\theta_\zs{j}}
\left(f(x,\theta)U(\theta)\right)\d \theta\,\nu_\zs{\xi}(\d x)\\[2mm]
&=\int_{\cX\times\bbr^{d-1}}\left(\int_{\bbr}\,
g^{\prime}_\zs{j}(\theta)\,
f(x,\theta)U(\theta)\d \theta_\zs{j}\right)\left(\prod_{i\neq j}\d \theta_i\right)\,\nu_\zs{\xi}(\d x)
=\eta_\zs{j}\,.
\end{align*}
Now by the Bouniakovskii-Cauchy-Schwarz inequality
we obtain the following lower bound for the quadratic risk
$$
\wt{\E}(\wh{g}_\zs{n}-g(\theta))^2\ge
\frac{\eta^2_\zs{j}}{\wt{\E}\wt{U}_\zs{j}^2}\,.
$$
To study the denominator in the left hand of this
inequality note that in view of the representation
 \eqref{sec:App.7}
$$
\frac{1}{f(y,\theta)}
\frac{\partial\,f(y,\theta)}{\partial\theta_\zs{j}}
=\frac{1}{\sigma_\zs{1}}\,
\int^{n}_\zs{0}\,\psi_\zs{j}(t)\,\d w_\zs{t}\,.
$$
Therefore, for each $\theta\in\bbr^d$,
$$
\E_\zs{\theta}\,
\frac{1}{f(y,\theta)}
\frac{\partial\,f(y,\theta)}{\partial\theta_\zs{j}}
\,
=0
$$
and
$$
\E_\zs{\theta}\,
\left(
\frac{1}{f(y,\theta)}
\frac{\partial\,f(y,\theta)}{\partial\theta_\zs{j}}
\right)^2
=\,
\frac{1}{\sigma^{2}_\zs{1}}
\int^{n}_\zs{0}\,\psi^2_\zs{j}(t)\d t
=
\frac{n}{\sigma^{2}_\zs{1}}
\Vert\psi\Vert^{2}
\,.
$$
Using equality
$$
\wt{U}_\zs{j}=
\frac{1}{f(x,\theta)}
\frac{\partial\,f(x,\theta)}{\partial\theta_\zs{j}}
+
\frac{1}{U(\theta)}
\frac{\partial\,U(\theta))}{\partial\theta_\zs{j}}
 \,,
$$
we get
$$
\wt{\E}\wt{U}_\zs{j}^2=
\frac{n}{\sigma^{2}_\zs{1}}\,\Vert\psi\Vert^{2}
+\,I_\zs{j}\,.
$$
Hence
Lemma~\ref{Le.sec:App.3}.
\endproof

\bigskip
\section{Proofs}\label{sec:Prf}

\subsection{Proof of Theorem \ref{Th.sec:Imp.1}}

Consider the quadratic error of the estimate \eqref{sec:Imp.11}
\begin{align*}
\|S^{*}_\zs{\lambda}-S\|^2&=\sum^{n}_\zs{j=1}(\lambda(j)\theta^{*}_\zs{\lambda,j}-\theta_\zs{j})^2=
\sum^{d}_\zs{j=1}(\lambda(j)\theta^{*}_\zs{\lambda,j}-\theta_\zs{j})^2+\sum^{n}_\zs{j=d+1}(\lambda(j)\wh{\theta}_\zs{j}-\theta_\zs{j})^2\\
&=\sum^{n}_\zs{j=1}(\lambda(j)\wh{\theta}_\zs{j}-\theta_\zs{j})^2
+
c^{2}_\zs{n}
-
2
c_\zs{n}
\sum^{d}_\zs{j=1}(\wh{\theta}_\zs{j}-\theta_\zs{j})\frac{\wh{\theta}_\zs{j}}{\|\wt{\theta}_\zs{n}\|_\zs{d}}
\\
&=\|\wh{S}_\zs{\lambda}-S\|^2+
c^{2}_\zs{n}
-
2
c_\zs{n}
\sum^{d}_\zs{j=1}(\wh{\theta}_\zs{j}-\theta_\zs{j})
\iota_\zs{j}(\wt{\theta}_\zs{n})
\,,
\end{align*}
where  $\iota_\zs{j}(x)=x_\zs{j}/\|x\|_\zs{d}$
for $x=(x_\zs{j})_\zs{1\le j\le d}\in\bbr^{d}$.
Therefore, we can represent the risk for the improved estimate $S^{*}_\zs{\lambda}$ as
$$
\cR_\zs{Q}(S^{*}_\zs{\lambda},S)=\cR_\zs{Q}(\wh{S}_\zs{\lambda},S)+
c^{2}_\zs{n}
-
2
c_\zs{n}
\,
\E_\zs{Q,S}\,\sum^{d}_\zs{j=1}(\wh{\theta}_\zs{j}-\theta_\zs{j})\,J_\zs{j}
\,,
$$
where $J_\zs{j}=\E(\iota_\zs{j}(\wt{\theta}_\zs{n})(\wh{\theta}_\zs{j}-\theta_j)|\cG_\zs{n})$.
Now, taking into account that the vector $\wt{\theta}_\zs{n}=(\wh{\theta}_\zs{j})_\zs{1\le j\le d}$
is the $\cG_\zs{n}$ conditionally gaussian vector in $\bbr^{d}$   with mean $\wt{\theta}=(\theta_\zs{j})_\zs{1\le j\le d}$
and covariance matrix $n^{-1}\G_\zs{n}$, we obtain
\begin{equation*}
J_\zs{j}
=\int_{\mathbb{R}^d}\,\iota_\zs{j}(x)(x-\theta_j)\p(x|\cG_\zs{n})
\d x\,.
\end{equation*}
Here $\p(x|\cG_\zs{n})$ is the conditional
distribution density of the vector $\wt{\theta}_\zs{n}$, i.e.
\begin{equation*}
\p(x|\cG_\zs{n})=\frac{1}{(2\pi)^{d/2}\sqrt{\det\G_\zs{n}}}
\exp\left(-\frac{(x-\theta)'\,\G^{-1}_\zs{n}(x-\theta)}{2}\right)
\,.
 \end{equation*}
Recall, that the $\prime$ denotes the transposition.
Changing the variables by
$u=\G^{-1/2}_\zs{n}(x-\theta)$,  one finds that
\begin{equation*}\label{sec:2.2c}
J_\zs{j}=\frac{1}{(2\pi)^{d/2}}\sum_{l=1}^{d}
\g_\zs{j,l}\int_{\mathbb{R}^{d}}\tilde{\iota}_\zs{j}(u)u_{l}\exp\left(-\frac{\|u\|^{2}_\zs{d}}
{2}\right) \d u\,,
 \end{equation*}
where
$\wt{\iota}_\zs{j}(u)=\iota_\zs{j}(\G^{1/2}_\zs{n} u+\theta)$
and $\g_\zs{ij}$ denotes the $(i,j)$-th\ element of  $\G^{1/2}_\zs{n}$.
Furthermore, integrating by parts, the integral $J_\zs{j}$ can be
rewritten as
\begin{equation*}
J_\zs{j}=\sum_{l=1}^{d}\sum_{k=1}^{d}\E\left(\g_\zs{jl}
\,\g_\zs{kl}\,
\frac{\partial \iota_\zs{j}}{\partial
u_k}(u)|_{u=\wt{\theta}_\zs{n}}|\cG_\zs{n}\right)
\,.
\end{equation*}
In view of the inequality $z^{\prime}Az\leq\lambda_{max}(A)\|z\|^2$
and Proposition \ref{sec:Imp.Prop_2_1}
we obtain that
\begin{align*}
\Delta_{Q}(S)&=
c^{2}_\zs{n}
-
2
c_\zs{n}
n^{-1}\E_\zs{Q,S}\,\left(\frac{\tr \G_\zs{n}}{\|\wt{\theta}_\zs{n}\|_\zs{d}}-
\frac{\wt{\theta}_\zs{n}^{\prime}\G_\zs{n}\wt{\theta}_\zs{n}}{\|\wt{\theta}_\zs{n}\|^3}\right)
\\[2mm]
&
\le c^{2}_\zs{n}
-
2
c_\zs{n}
\,(d-1)\underline{\sigma}_n n^{-1}\E_\zs{Q,S}\,\frac{1}{\|\wt{\theta}_\zs{n}\|_\zs{d}}
\,.
\end{align*}
Moreover, using the Jensen inequality we can estimate the last expectation from below as
$$
\E_\zs{Q,S}\,(\|\wt{\theta}_\zs{n}\|_\zs{d})^{-1}=\E_\zs{Q,S}\,(\|\wt{\theta}+n^{-1/2}\wt{\xi}_\zs{n}\|_\zs{d})^{-1}
\geq
\,(\|\theta\|_\zs{d}+n^{-1/2}\E_\zs{Q,S}\|\wt{\xi}_\zs{n}\|_\zs{d})^{-1}
\,.
$$
From Proposition \ref{Pr.sec:Stc.1} and
 the condition
\eqref{sec:Ex.01-1}
we  obtain
$$
\E_\zs{Q,S}\|\wt{\xi}_\zs{n}\|^{2}_\zs{d}\le \overline{\sigma}_\zs{n}\,d\,.
$$
So, for $\Vert S\Vert\le r_\zs{n}$
$$
\E_\zs{Q,S}\,\|\wt{\theta}_\zs{n}\|^{-1}\geq
\left(r_\zs{n}+\sqrt{d\overline{\sigma}_\zs{n}/n}\right)^{-1}
$$
and, therefore,
$$
\Delta_{Q}(S)
\le c^{2}_\zs{n}
-
2
c_\zs{n}
\frac{(d-1)\underline{\sigma}_n}{\left(r_\zs{n}+\sqrt{d\overline{\sigma}_\zs{n}/n}\right)\,n}
=-c^{2}_\zs{n}
\,.
$$
Hence
Theorem \ref{Th.sec:Imp.1}.
\endproof

\subsection{Proof of Theorem \ref{sec:Mo.Th.1}}

\noindent
Using the definitions \eqref{sec:Mo.1}, \eqref{sec:Mo.2vartheta} and
 \eqref{sec:Mo.4}, we obtain that  for any $\lambda\in\Lambda$
\begin{align}\nonumber
\Er_\zs{n}(\lambda)\,&=\,J_\zs{n}(\lambda)+
2\,\sum^{n}_\zs{j=1}\,\lambda(j)\left(\theta^{*}_\zs{\lambda,j}\,\wh{\theta}_\zs{j}-\frac{\wh{\sigma}_\zs{n}}{n}
-\theta^{*}_\zs{\lambda,j}\,\theta_\zs{j}\right)\\[2mm]
\label{sec:Mo.11}
&+\,
\|S\|^2-\delta\wh{P}_\zs{n}(\lambda)\,.
\end{align}
Now we set
\begin{align}\nonumber
B_\zs{1,n}&(\lambda)=\frac{1}{\sqrt{n}}\sum^{n}_\zs{j=1}\,\lambda(j)g_\zs{\lambda}(j)\wh{\theta}_\zs{j}\xi_\zs{j}\,,\qquad
B_\zs{2,n}(\lambda)=\sum^{n}_\zs{j=1}\,\lambda(j)
\,\wt{\xi}_\zs{j}\\[2mm]\label{sec:Mo_Proof_B_M}
&\mbox{and}\qquad
M(\lambda)=\frac{1}{\sqrt{n}}\sum^{n}_\zs{j=1}\,\lambda(j)\theta_\zs{j}\xi_\zs{j}
\,,
\end{align}
where
$g_\zs{\lambda}(j)=(c_\zs{n}(\lambda)/\vert\wt{\theta}\vert_\zs{d})\,\Chi_\zs{\{1\le j\le d\}}$ and $\wt{\xi}_\zs{j}=\xi^2_\zs{j}-\E_\zs{Q} \xi^2_\zs{j}$.
Taking into account the definition \eqref{sec:Mo.5_whP}, we can rewrite \eqref{sec:Mo.11} as
\begin{align}\nonumber
\Er_\zs{n}(\lambda)&=J_\zs{n}(\lambda)+2\frac{\sigma-\wh{\sigma}_\zs{n}}{n}L(\lambda)+
2\,M(\lambda)
\\[2mm]\label{sec:Mo.12}
&-2B_\zs{1,n}(\lambda)\,+2\sqrt{P_\zs{n}(\lambda)}\frac{B_\zs{2,n}(\overline{\lambda})}{\sqrt{\sigma n}}+\,
\|S\|^2-\delta\wh{P}_\zs{n}(\lambda)\,,
\end{align}
where
 the function $L(\lambda)$ is defined in
\eqref{sec:Imp.6+0},
 $\overline{\lambda}=\lambda/|\lambda |_\zs{n}$.
 Let $\lambda_0=(\lambda_0(j))_{1\le n}$ be a fixed sequence in $\Lambda$ and $\lambda^*$ be as in \eqref{sec:Mo.6}.
Substituting $\lambda_0$ and $\lambda^*$ in \eqref{sec:Mo.12}, we consider the difference
\begin{align*}
\Er_\zs{n}(\lambda^*)-\Er_\zs{n}(\lambda_0)&\leq 2\frac{\sigma-\wh{\sigma}_\zs{n}}{n}L(\varpi)+2M(\varpi)-2B_\zs{1,n}(\lambda^*)+2B_\zs{1,n}(\lambda_0)\\[2mm]
&+2\sqrt{P_\zs{n}(\lambda^*)}\frac{B_\zs{2,n}(\overline{\lambda^*})}{\sqrt{\sigma n}}
-2\sqrt{P_\zs{n}(\lambda_0)}\frac{B_\zs{2,n}(\overline{\lambda_0})}{\sqrt{\sigma n}}\\[2mm]
&
-\delta\wh{P}_\zs{n}(\lambda^*)
+\delta\wh{P}_\zs{n}(\lambda_0)\,,
\end{align*}
where $\varpi=\lambda^*-\lambda_0\in\Lambda_\zs{1}$ and
\begin{equation}
\label{setLambda-1}
\Lambda_\zs{1}=\Lambda-\lambda_\zs{0}=\left\{
\lambda-\lambda_\zs{0}\,,\,\lambda\in\Lambda
\right\}
\,.
\end{equation}
Note that $|L(\varpi)|\leq 2\vert \Lambda\vert_\zs{*}$.
Moreover, note also that
\begin{equation}\label{sec:Mo.13_Ub++c-n}
\sum^{n}_\zs{j=1}\,g^{2}_\zs{\lambda}(j)\,
\wh{\theta}^{2}_\zs{j}
=c^{2}_\zs{n}(\lambda)
\le \frac{c^{*}_\zs{n}}{n}
\,,
\end{equation}
where $c^{*}_\zs{n}$ is defined
in \eqref{sec:Mo.9+1_psi}.
Therefore,
 through the  Cauchy--Schwarz inequality we can estimate the term $B_\zs{1,n}(\lambda)$ as
$$
|B_\zs{1,n}(\lambda)|\le
\frac{|\lambda |_\zs{n}}{\sqrt{n}}c_\zs{n}(\lambda)
\left(\sum_\zs{j=1}^{n}
\overline{\lambda}^{2}(j)
\,
\xi^{2}_\zs{j}
\right)^{1/2}
=
\frac{|\lambda |_\zs{n}}{\sqrt{n}}c_\zs{n}(\lambda)
\left(
\sigma+ B_\zs{2,n}(\overline{\lambda}^{2})
\right)^{1/2}
\,,
$$
where $x^{2}=(x^{2}(j))_\zs{1\le j\le n}$ for $x\in\bbr^{n}$.
So, applying the elementary inequality
\begin{equation}\label{sec:Mo.13eleq-0}
2|ab|\leq \varepsilon a^2+\varepsilon^{-1} b^2
\end{equation}
with some arbitrary $\varepsilon>0$, we get
$$
2|B_\zs{1,n}(\lambda)|
\leq \varepsilon P_\zs{n}(\lambda)+\frac{c^{*}_\zs{n}}{\varepsilon\sigma n}
(\sigma+B^{*}_\zs{2})\,.
$$
Moreover, by the same method we estimate the term $B_\zs{2,n}$, i.e.
$$
2\sqrt{P_\zs{n}(\lambda)}\frac{B_\zs{2,n}(\overline{\lambda})}{\sqrt{\sigma n}}
\leq \varepsilon P_\zs{n}(\lambda)+\frac{B_\zs{2,n}^2(\overline{\lambda})}{\varepsilon\sigma n}
\leq \varepsilon P_\zs{n}(\lambda)+\frac{B^{*}_\zs{2}}{\varepsilon \sigma n}\,,
$$
where
$$
B^{*}_\zs{2}
=
\max_\zs{\lambda\in\Lambda}\,
\left(
B_\zs{2,n}^2(\overline{\lambda})
+
B_\zs{2,n}^2(\overline{\lambda}^{2})
\right)
\,.
$$
 Note that from  Proposition \ref{Pr.sec:Stc.4}
we obtain that
\begin{equation}\label{sec:Mo.13_Ub}
\E_\zs{Q}\,B^{*}_\zs{2}
\le
\sum_\zs{\lambda\in\Lambda}\,
\left(
\E_\zs{Q}B_\zs{2,n}^2(\overline{\lambda})
+
\E_\zs{Q}B_\zs{2,n}^2(\overline{\lambda}^{2})
\right)
\le
2(2+\overline{\phi}_\zs{n}^4\Pi(x^4))
\sigma^2\nu_n\,.
\end{equation}
Using the bounds above, one has
\begin{multline*}
\Er_\zs{n}(\lambda^*)\leq\Er_\zs{n}(\lambda_0)+\frac{4\vert\Lambda\vert_\zs{n} |\wh{\sigma}_\zs{n}-\sigma|}{n}
+2M(\varpi)
\\[2mm]
+\frac{2}{\varepsilon}\,\frac{c^{*}}{n\sigma}(\sigma+B^{*}_\zs{2})
+\frac{2}{\varepsilon}\,
\frac{B^{*}_\zs{2}}{n\sigma}
\\[2mm]
+2 \varepsilon P_\zs{n}(\lambda^*)
+2\varepsilon P_\zs{n}(\lambda_0)
-\delta\wh{P}_\zs{n}(\lambda^*)+\delta\wh{P}_\zs{n}(\lambda_0)\,.
\end{multline*}
The setting $\varepsilon=\delta/4$
and the estimating where this is possible $\delta$ by $1$
 in this inequality
imply
\begin{multline*}
\Er_\zs{n}(\lambda^*)\leq\Er_\zs{n}(\lambda_0)+
 \frac{5\vert\Lambda\vert_\zs{n}|\wh{\sigma}_\zs{n}-\sigma|}{n}
+2M(\varpi)
\\[2mm]
+\frac{16 (c^{*}_\zs{n}+1)(\sigma+B^{*}_\zs{2})}{\delta n\sigma}
-\frac{\delta}{2}\wh{P}_\zs{n}(\lambda^*)+\frac{\delta}{2} P_\zs{n}(\lambda_0)
+\delta\wh{P}_\zs{n}(\lambda_0)\,.
\end{multline*}
Moreover, taking into account here that
$$
\vert
\wh{P}_\zs{n}(\lambda_0)
-
P_\zs{n}(\lambda_0)
\vert
\le
\frac{\vert\Lambda\vert_\zs{n}|\wh{\sigma}_\zs{n}-\sigma|}{n}
$$
and that $\delta<1/2$,
we obtain that
\begin{multline}\label{sec:Mo.14}
\Er_\zs{n}(\lambda^*)\leq\Er_\zs{n}(\lambda_0)+
 \frac{6\vert\Lambda\vert_\zs{n}|\wh{\sigma}_\zs{n}-\sigma|}{n}
+2M(\varpi)
\\[2mm]
+\frac{16 (c^{*}_\zs{n}+1)(\sigma+B^{*}_\zs{2})}{\delta n\sigma}
-\frac{\delta}{2} P_\zs{n}(\lambda^*)+\frac{3\delta}{2} P_\zs{n}(\lambda_0)\,.
\end{multline}

Now we examine the third term in the right-hand side of this inequality. Firstly we note that
\begin{equation}
\label{upper_bound_M-+01}
2|M(\varpi)|\leq\varepsilon\|S_\zs{\varpi}\|^2+\frac{Z^*}{n\varepsilon}\,,
\end{equation}
where $S_\zs{\varpi}=\sum^{n}_\zs{j=1}\,\varpi_\zs{j}\theta_\zs{j}\phi_\zs{j}$
and
$$
Z^*=\sup_{x\in\Lambda_1}\frac{nM^2(x)}{\|S_x\|^2}
\,.
$$
We remind that the set $\Lambda_\zs{1}$ is defined in
\eqref{setLambda-1}.
Using Proposition
\ref{Pr.sec:Stc.1}
we can obtain that for any fixed $x=(x_\zs{j})_\zs{1\le j\le n}\in\bbr^{n}$
\begin{equation}
\label{M^2+11-00}
\E\,M^2(x)=\frac{\E\,I^{2}_\zs{n}\left(S_\zs{x}\right)}{n^{2}}
=\frac{\sigma \Vert S_\zs{x}\Vert^{2}}{n}
=\frac{\sigma}{n}\,
\sum^{n}_\zs{j=1}\,x^{2}_\zs{j}\,\theta^{2}_\zs{j}
\end{equation}
and, therefore,
\begin{equation}
\label{up-Z*-00}
\E_\zs{Q}Z^*
\le
\sum_{x\in\Lambda_1}\frac{nM^2(x)}{\|S_x\|^2}
\leq \sigma
\nu_n
\,.
\end{equation}
Moreover,  the norm $\Vert S^{*}_\zs{\lambda^{*}}-S^{*}_\zs{\lambda_\zs{0}}\Vert$ can be estimated from below as
\begin{align*}
\Vert S^{*}_\zs{\lambda}-S^{*}_\zs{\lambda_\zs{0}}\Vert^{2}
&=
\sum^{n}_\zs{j=1}
(\varpi(j)+\beta(j))^{2}\wh{\theta}^{2}_\zs{j}
\\[2mm]
&\ge
\|\wh{S}_\zs{\varpi}\|^2
+2
\sum^{n}_\zs{j=1}
\varpi(j)\beta(j)\wh{\theta}^{2}_\zs{j}\,,
\end{align*}
where $\beta(j)=\lambda_0(j)g_j(\lambda_0)-\lambda(j)g_j(\lambda)$.
Therefore, in view of \eqref{sec:Imp.4}
\begin{align*}
\|S_\zs{\varpi}\|^2&-\|S^{*}_\zs{\lambda}-S^{*}_\zs{\lambda_\zs{0}}\|^2
\le
\|S_\zs{\varpi}\|^2-\|\wh{S}_\zs{\varpi}\|^2
-2\sum^{n}_\zs{j=1}\,\varpi(j)\beta(j)
\wh{\theta}_\zs{j}^2
\\[2mm]
&
\le
-2M(\varpi^{2})-2\sum^{n}_\zs{j=1}\,\varpi(j)\beta(j)
\wh{\theta}_\zs{j}\theta_j-
\frac{2}{\sqrt{n}}\Upsilon(\varpi)
\,,
\end{align*}
where $\Upsilon(\lambda)=\sum^{n}_\zs{j=1}\,\lambda(j)\beta(j)
\wh{\theta}_\zs{j}\xi_\zs{j}$. Note that the first term in this inequality we can estimate as
$$
2M(\varpi^{2})\le \varepsilon\|S_\zs{\varpi}\|^2+\frac{Z_1^*}{n\varepsilon}
\quad\mbox{and}\quad
Z^*_1=\sup_{x\in\Lambda_1}\frac{n M^2(x^{2})}{\|S_x\|^2}
\,.
$$
Note that, similarly to \eqref{up-Z*-00} we can estimate the last term as
$$
\E_\zs{Q}Z_1^*\leq \sigma\nu_n\,.
$$
From this it follows that for any $0<\varepsilon<1$
\begin{align}
\nonumber
\|S_\zs{\varpi}\|^2
\le
\frac{1}{1-\varepsilon}
&\left(
\|S^{*}_\zs{\lambda}-S^{*}_\zs{\lambda_\zs{0}}\|^2
+\frac{Z_1^*}{n\varepsilon}\right.
 \\[2mm] \label{upper-bound-000}
&
\left.
-2\sum^{n}_\zs{j=1}\,\varpi(j)\beta(j)
\wh{\theta}_\zs{j}\theta_j-
\frac{2\Upsilon(\varpi)}{\sqrt{n}}
\right)
\,.
\end{align}
Moreover, note now that the property
 \eqref{sec:Mo.13_Ub++c-n} yields
\begin{equation}
\label{theta-whg-upperb-00}
\sum^{n}_\zs{j=1}\,
\beta^{2}(j)
\wh{\theta}^{2}_\zs{j}
\le
2
\sum^{n}_\zs{j=1}\,g^{2}_\zs{\lambda}(j)\,
\wh{\theta}^{2}_\zs{j}
+
2\sum^{n}_\zs{j=1}\,g^{2}_\zs{\lambda_\zs{0}}(j)\,
\wh{\theta}^{2}_\zs{j}
\le
\frac{4c^{*}}{\varepsilon n}\,.
\end{equation}
Taking into account that $\vert\varpi(j)\vert\le 1$ and using the inequality \eqref{sec:Mo.13eleq-0}, we get that
 for any $\varepsilon>0$
$$
2\left|\sum^{n}_\zs{j=1}\,\varpi(j)
\beta(j)
\wh{\theta}_\zs{j}\theta_j\right|\leq\varepsilon\|S_\zs{\varpi}\|^2
+\frac{4c^{*}}{\varepsilon n}
\,.
$$
To estimate the last term in the right hand of  \eqref{upper-bound-000} we use first
the Cauchy -- Schwarz inequality
and then the bound
\eqref{theta-whg-upperb-00}, i.e.
\begin{align*}
\frac{2}{\sqrt{n}}\vert \Upsilon(\lambda)\vert
&\le
\frac{2\vert\lambda\vert_\zs{n}}{\sqrt{n}}\left(\sum^{n}_\zs{j=1}\,\beta^{2}(j)
\wh{\theta}^{2}_\zs{j}\right)^{1/2}
\left(
\sum^{n}_\zs{j=1}\,\bar{\lambda}^{2}(j)\,\xi^{2}_\zs{j} \right)^{1/2}\\[2mm]
&
\le
\varepsilon P_\zs{n}(\lambda)
+
\frac{c^{*}}{n\varepsilon\sigma}
\sum^{n}_\zs{j=1}\,\bar{\lambda}^{2}(j)\,\xi^{2}_\zs{j}
\le
\varepsilon P_\zs{n}(\lambda)
+
\frac{c^{*}(\sigma+B^{*}_\zs{2})}{n\varepsilon\sigma}
\,.
\end{align*}
Therefore,
\begin{align*}
\frac{2}{\sqrt{n}}\vert \Upsilon(\varpi)\vert
&\le
\frac{2}{\sqrt{n}}\vert \Upsilon(\lambda^{*})\vert
+\frac{2}{\sqrt{n}}\vert \Upsilon(\lambda_\zs{0})\vert\\[2mm]
&\le
\varepsilon P_\zs{n}(\lambda^{*})
+\varepsilon P_\zs{n}(\lambda_\zs{0})
+
\frac{2c^{*}(\sigma+B^{*}_\zs{2})}{n\varepsilon\sigma}
\,.
\end{align*}
So, using all these bounds in \eqref{upper-bound-000},
we obtain that
\begin{align*}
\|S_\zs{\varpi}\|^{2}
&\le
\frac{1}{(1-\varepsilon)}\biggl(\frac{Z_1^*}{n\varepsilon}+
\|S_\zs{\lambda^*}^*-S_\zs{\lambda_0}^*\|^2
+\frac{6c^{*}_\zs{n}(\sigma+B^{*}_\zs{2})}{n\sigma\varepsilon}
\\[2mm]
&+\varepsilon P_\zs{n}(\lambda^*)+\varepsilon P_\zs{n}(\lambda_0)\biggr)\,.
\end{align*}
Using in the inequality \eqref{upper_bound_M-+01} this bound and the estimate
$$
\|S_\zs{\lambda^*}^*-S_\zs{\lambda_0}^*\|^2\leq
2(\Er_\zs{n}(\lambda^*)+\Er_\zs{n}(\lambda_0))\,,
$$
we obtain
\begin{align*}
2|M(\varpi)|&\le
\frac{Z^*+Z_1^*}{n(1-\varepsilon)\varepsilon}
+
\frac{2\varepsilon(\Er_\zs{n}(\lambda^*)+\Er_\zs{n}(\lambda_0))}{(1-\varepsilon)}
\\[2mm]
&
+
\frac{6c^{*}_\zs{n}(\sigma+B^{*}_\zs{2})}{n\sigma(1-\varepsilon)}
+\frac{\varepsilon^{2}}{1-\varepsilon}\left( P_\zs{n}(\lambda^*)+ P_\zs{n}(\lambda_0)\right)
\,.
\end{align*}
Choosing here $\varepsilon\le \delta/2<1/2$ we obtain that
\begin{align*}
2|M(\varpi)|&\le
\frac{2(Z^*+Z_1^*)}{n\varepsilon}
+
\frac{2\varepsilon(\Er_\zs{n}(\lambda^*)+\Er_\zs{n}(\lambda_0))}{(1-\varepsilon)}
\\[2mm]
&
+
\frac{12c^{*}_\zs{n}(\sigma+B^{*}_\zs{2})}{n\sigma}
+
\varepsilon\left( P_\zs{n}(\lambda^*)+ P_\zs{n}(\lambda_0)\right)
\,.
\end{align*}
From here and \eqref{sec:Mo.14}, it follows that
\begin{align*}
\Er_\zs{n}(\lambda^*) &\leq\frac{1+\varepsilon}{1-3\varepsilon}\Er_\zs{n}(\lambda_0)
+ \frac{6\vert\Lambda\vert_\zs{n}|\wh{\sigma}_\zs{n}-\sigma|}{n(1-3\varepsilon)}
\\[2mm]
&
+\frac{28(1+c^{*}_\zs{n})(B^{*}_\zs{2}+\sigma)}{\delta(1-3\varepsilon)n\sigma}
+\frac{2(Z^*+Z_1^*)}{n(1-3\varepsilon)}
+\frac{2\delta P_\zs{n}(\lambda_0)}{1-3\varepsilon}.
\end{align*}
Choosing here $\varepsilon=\delta/3$
and estimating $(1-\delta)^{-1}$ by $2$ where this is possible,
 we get
\begin{align*}
\Er_\zs{n}(\lambda^*) &\leq\frac{1+\delta/3}{1-\delta}\Er_\zs{n}(\lambda_0)
+ \frac{12\vert\Lambda\vert_\zs{n}|\wh{\sigma}_\zs{n}-\sigma|}{n}
\\[2mm]
&
+\frac{56(1+c^{*}_\zs{n})(B^{*}_\zs{2}+\sigma)}{\delta n\sigma}
+\frac{4(Z^*+Z_1^*)}{n}
+\frac{2\delta P_\zs{n}(\lambda_0)}{1-\delta}\,.
\end{align*}
Taking the expectation and using the upper bound for $P_\zs{n}(\lambda_0)$  in Lemma~\ref{Lem.A.1} with $\varepsilon=\delta$ yields
$$
\cR_\zs{Q}(S^*,S)\leq\frac{1+5\delta}{1-\delta}\cR_\zs{Q}(S^*_\zs{\lambda_0},S)
+\frac{\check{\U}_\zs{Q,n}}{n\delta}
+
 \frac{12\vert\Lambda\vert_\zs{n}\E_\zs{Q}|\wh{\sigma}_\zs{n}-\sigma|}{n}
\,,
$$
where
$
\check{\U}_\zs{Q,n}=
56(1+c^{*}_\zs{n})
(2(6+\overline{\phi}_\zs{n}^4\Pi(x^4))
\sigma\nu_n
+
1)
+2
c^{*}_\zs{n}$.
The inequality holds for each $\lambda_0\in\Lambda$, this  implies  Theorem \ref{sec:Mo.Th.1}. \endproof

\bigskip

\subsection{Proof of Theorem \ref{Th.sec: Ae.2}}

Firstly, note, that for any fixed $Q\in \cQ_\zs{n}$
\begin{equation}\label{sec:Lo.1-0}
\sup_\zs{S\in W_\zs{k,r}}\,\cR^{*}_\zs{n}(\wh{S}_\zs{n},S)
\ge
\,
\sup_\zs{S\in W_\zs{k,r}}\,
\cR_\zs{Q}(\wh{S}_\zs{n},S)
\,.
\end{equation}
Now for any fixed $0<\ve<1$ we set
\begin{equation}\label{sec:Lo.1}
d=d_\zs{n}=\left[\frac{k+1}{k}v^{1/(2k+1)}_\zs{n}l_\zs{k}(r_\zs{\ve})\right]
\quad\mbox{and}\quad
r_\zs{\ve}=(1-\ve) r\,.
\end{equation}
Next we approximate
the unknown function by a trigonometric series with $d=d_\zs{n}$ terms, i.e.
for any array $z=(z_\zs{j})_\zs{1\le j\le d_\zs{n}}$,
 we set
\begin{equation*}\label{sec:Lo.4}
S_\zs{z}(x)=\sum_{j=1}^{d_\zs{n}}\,z_\zs{j}\,\phi_\zs{j}(x)
\,.
\end{equation*}

\noindent
To define the Bayesian risk we  choose a prior distribution on $\bbr^{d}$
as
\begin{equation*}\label{sec:Lo.5}
\kappa=(\kappa_\zs{j})_\zs{1\le j\le d_\zs{n}}
\quad\mbox{and}\quad
\kappa_\zs{j}=s_\zs{j}\,\eta_\zs{j}\,,
\end{equation*}
where $\eta_\zs{j}$ are i.i.d. gaussian $\cN(0,1)$ random variables and
the coefficients
$$
s_\zs{j}=\sqrt{\frac{s^*_\zs{j}}{v_\zs{n}}}
\quad\mbox{and}\quad
s^{*}_\zs{j}\,
=
\left( \frac{d_\zs{n}}{j}
\right)^{k}
-
1
\,.
$$
Furthermore, for any function $f$, we denote by $\p(f)$ its projection
in $\L_\zs{2}[0,1]$
 onto
 $W_\zs{k,r}$, i.e.
 $$
\p(f)=\hbox{\rm Pr}_\zs{W_\zs{k,r}}(f)\,.
$$
Since $W_\zs{k,r}$ is a convex set, we obtain
$$
\|\wh{S}-S\|^2\ge\|\wh{\p}-S\|^2
\quad\mbox{with}\quad
\wh{\p}=\p(\wh{S})
\,.
$$
Therefore,
$$
\sup_\zs{S\in W_\zs{k,r}}\,\cR(\wh{S},S)
\ge\,
\int_{\{z\in\bbr^d\,:\,S_\zs{z}\in W_\zs{k,r}\}}\,
\E_\zs{S_\zs{z}}\|\wh{\p}-S_\zs{z}\|^2\,\mu_{\kappa}(\d z)
\,.
$$
Using the distribution $\mu_\zs{\kappa}$ we introduce
 the following Bayes risk
$$
\wt{\cR}_\zs{Q}(\wh{S})=
\int_\zs{\bbr^d}\cR_\zs{Q}(\wh{S},S_\zs{z})\,
\mu_\zs{\kappa}(\d z)\,.
$$
Taking into account now that $\|\wh{\p}\|^2\le r$
 we obtain
\begin{equation}\label{sec:Lo.12}
\sup_\zs{S\in W_\zs{k,r}}\,
\cR_\zs{Q}(\wh{S},S)
\,
\ge\,
\wt{\cR}_\zs{Q}(\wh{\p})
-2\,
\R_\zs{0,n}
\end{equation}
with
$$
\R_\zs{0,n}=
\int_\zs{ \{z\in\bbr^d\,:\,S_\zs{z}\notin W_\zs{k,r}\}}\,
\,
(r+\|S_\zs{z}\|^2)\,
\mu_\zs{\kappa}(\d z)
\,.
$$
Therefore, in view of \eqref{sec:Lo.1-0}
\begin{equation}\label{sec:Lo.1+10}
\sup_\zs{S\in W_\zs{k,r}}\,\cR^{*}_\zs{n}(\wh{S}_\zs{n},S)
\ge
\,
\sup_\zs{Q\in\cQ_\zs{n}}\,
\wt{\cR}_\zs{Q}(\wh{\p})
-2\,
\R_\zs{0,n}
\,.
\end{equation}
In
Lemma \ref{Le.sec:App.3+1}
we studied the last term in this inequality.
Now it is easy to see that
$$
\|\wh{\p}-S_z\|^2 \ge
\sum_{j=1}^{d_\zs{n}}\,
(\wh{z}_\zs{j}-z_\zs{j})^2
\,,
$$
where $\wh{z}_\zs{j}=\int^{1}_\zs{0}\,\wh{\p}(t)\,\phi_\zs{j}(t)\d t.$ So, in view of Lemma~\ref{Le.sec:App.3}
and reminding
that
$v_\zs{n}=n/\overline{\sigma}_\zs{n}$
 we obtain
\begin{align*}
\sup_\zs{Q\in\cQ_\zs{n}}
\wt{\cR}_\zs{Q}(\wh{\p})\,
&\ge\,
\sup_\zs{0<\sigma^{2}_\zs{1}\le \sigma^{*}}
\sum_{j=1}^{d_\zs{n}}\,\frac{1}
{n\sigma^{-2}_\zs{1}+v_\zs{n}\,(s^{*}_\zs{j})^{-1}}
\\[2mm]
&
=
\frac{1}{v_\zs{n}}\,
\sum_{j=1}^{d_\zs{n}}\,\frac{s^{*}_\zs{j}}
{s^{*}_\zs{j}+\,1}
=
\frac{1}{v_\zs{n}}\,
\sum_{j=1}^{d_\zs{n}}\,
\left(
1
-
\frac{j^k}{d^k_\zs{n}}
\right)
\,.
\end{align*}
Therefore, using now the definition \eqref{sec:Lo.1},
Lemma \ref{Le.sec:App.3+1} and the inequality
\eqref{sec:Lo.1+10}
we obtain that
$$
\liminf_\zs{n\to\infty}\inf_\zs{\wh{S}\in\Pi_\zs{n}}\,v^{\frac{2k}{2k+1}}_\zs{n}\,
\sup_\zs{S\in W_\zs{k,r}}\,\cR^{*}_\zs{n}(\wh{S}_\zs{n},S)
\ge\,
(1-\ve)^{\frac{1}{2k+1}}\,
l_\zs{k}(r_\zs{\ve})\,.
$$
Taking here limit as $\ve\to 0$ we come to the Theorem~\ref{Th.sec: Ae.2}.
\endproof

\bigskip

\subsection{Proof of Theorem \ref{Th.sec: Ae.1}}

This theorem follows from Theorem \ref{Th.sec:Imp.1} and Theorem 3.1 in \cite{KonevPergamenshchikov2009b}.
\endproof

\bigskip

\section{Conclusion}

In the conclusion we would like to emphasize that in this paper    we develop a  new
model selection method based on the improved  versions of the least squares  estimates.
It turns out that the improvement effect in the nonparametric estimation
given in
\eqref{sec:Imp.11+1}
 is more  important than  for the parameter estimation problems
since the  accuracy improvement is proportional to the parameter dimension $d$ which goes to infinity for nonparametric models.
Recall that,
 the improved estimation methods was
 usually used
 for the parametric estimation problem only, where the parameter dimension $d$ is always fixed  (see, for example,  \cite{FourdrinierStrawderman1996}).
 Therefore,  the benefit in the non-asymptotic quadratic accuracy
 from the application of the improved estimation methods is more significant in  statistical nonparametric signal processing.
Moreover, for the proposed improved model selection procedures we obtain the  sharp oracle inequalities. It should be emphasized that
in this paper we obtain these inequalities without conditions on the jumps, i.e. without assumption that the L\'evy measure is finite.
To this end we developed a special analytical tool in Proposition \ref{Pr.sec:Stc.3} to study the non-asymptotic properties for the
corresponding
stochastic integrals with respect to the process  \eqref{sec:In.1+1}.
Moreover, asymptotically, as $n$ goes to infinity,
we shown the adaptive
efficiency
for
the  improved model selection procedures.
This is the meaning, that the proposed shrinkage  model selection procedures
have the benefit with respect to the least squares estimator
in the non-asymptotic accuracy and asymptotically they
possess the same efficient properties as the least squares methods.
Moreover, the behavior of the constructed procedures
 is illustrated by the numerical simulations in Section \ref{sec:Sim}.

\bigskip

{\bf Acknowledgements.}
This work was partially supported  by the research project no. 2.3208.2017/4.6
(the  Ministry of Education and Science of the Russian Federation), RFBR Grant 16-01-00121 A and
"The Tomsk State University competitiveness improvement programme"\ Grant 8.1.18.2018.
The work of the last author was  partially supported
 by
  the Russian Federal Professor program (project no. 1.472.2016/1.4, the  Ministry of Education and Science of the Russian Federation)
  and by the European research project XterM - Feder, University of Havre (France).

\bigskip

\medskip



\bigskip

\renewcommand{\theequation}{A.\arabic{equation}}
\renewcommand{\thetheorem}{A.\arabic{theorem}}
\renewcommand{\thesubsection}{A.\arabic{subsection}}
\section{Appendix}\label{sec:A}
\setcounter{equation}{0}
\setcounter{theorem}{0}

\bigskip

\subsection{Proof of Proposition \ref{sec:Imp.Prop_2_1}}

Using \eqref{sec:In.1+1} we put for any square integrated  functions $f$
\begin{equation*}\label{sec:Stc.7}
I_\zs{t}^{(1)}(f)=\int_0^t f(s) \d w_s
\quad\mbox{and}\quad
I_\zs{t}^{(2)}(f)=\int_0^t f(s) \d z_s\,.
\end{equation*}
From here and \eqref{sec:Imp.4} we can see that the vector $\tilde{\xi}$ has the conditionally Gaussian distribution with respect to $\cG_n$ with zero mean and its covariance matrix $\G_\zs{n}$ can be rewritten as
$$
\G_\zs{n}=\sigma_1^2\,I_\zs{d}+\sigma_2^2\cD_\zs{n},
$$
where $I_\zs{d}$ is the identity matrix and the $(i,j)$-th element of the matrix $\cD_\zs{n}$ is defined as
$\E(I_n^{(2)}(\phi_\zs{i})I_n^{(2)}(\phi_\zs{j})|\cG_n)$. Using the celebrated inequality
of Lidskii and Wieland (see, for example, in
\cite{MarchallOlkin1979}, G.3.a., p.334)
 we obtain
\begin{equation*}\label{eq1.24}
\tr\G_n-\lambda_{\max}(\G_n)\ge\sigma_1^{2}(
\tr\,I_\zs{d}-\lambda_{\max}(I_\zs{d}))
\quad\mbox{a.s.}
\end{equation*}
Now,  using \eqref{sec:Ex.01-1} we come to desire results.
\endproof

\bigskip

\subsection{Proof of Proposition \ref{sec:Mo.Prop.1}}

We use here the same method as in \cite{KonevPergamenshchikov2009a}.
Using the equality \eqref{sec:Imp.4} for the trigonometric basis, we get
\begin{equation*}\label{sec:Mo.Prop_1.1}
\wh{t}_\zs{j}=
t_\zs{j}+
\frac{1}{\sqrt{n}}
\xi_\zs{j}\,,
\end{equation*}
where
$$
t_\zs{j}=
\int^{1}_\zs{0}\,S(u)\,\Trg_\zs{j}(u)\d u
\quad\mbox{and}\quad
\xi_\zs{j}=
\frac{1}{\sqrt{n}}
\int^{n}_\zs{0}\,\Trg_\zs{j}(u)\,\d \xi_\zs{u}
\,.
$$
Therefore, the estimator \eqref{sec:Mo.3} can be represented as
\begin{equation}\label{sec:Mo.Prop_010-04}
\wh{\sigma}_\zs{n}=\sum_\zs{j=[\sqrt{n}]+1}^n\,t^{2}_\zs{j}
+
2\frac{1}{\sqrt{n}}
M_\zs{n}
+
\frac{1}{n}
\sum^n_\zs{j=[\sqrt{n}]+1}\,\xi^{2}_\zs{j}
\,,
\end{equation}
where $M_\zs{n}=\sum^n_\zs{j=[\sqrt{n}]+1}\,t_\zs{j}\,\xi_\zs{j}$.
Note that for the continuously
differentiable functions (see, for example, Lemma A.6 in \cite{KonevPergamenshchikov2009a})
the Fourrier coefficients $(t_\zs{j})$
for any $m\ge 1$ satisfy  the following inequality
\begin{equation*}\label{sec:Mo.2-1-04}
\sum^{\infty}_\zs{j=m+1}\,t^2_\zs{j}
\le
\frac{4}{m}
\left(\int^{1}_\zs{0}\vert\dot{S}(t)\vert \d t\right)^{2}
\le
\frac{4}{m}
\Vert\dot{S}\Vert^{2}
\end{equation*}
and $\dot{S}=\d S/\d t$.  The second term in \eqref{sec:Mo.Prop_010-04}
can be estimated through
 the equality \eqref{M^2+11-00}, i.e.
$$
\E_\zs{Q}\,M^{2}_\zs{n}=\frac{\sigma}{n}
 \sum^{n}_\zs{j=[\sqrt{n}]+1}\,t^{2}_\zs{j}
\le
\frac{4\sigma}{n\sqrt{n}}\,
\Vert\dot{S}\Vert^{2}\,.
$$
Moreover, taking into account that
 the expectation $\E_\zs{Q}\,\xi^{2}_\zs{j}=\sigma$
we can represent the last term in \eqref{sec:Mo.Prop_010-04} as
$$
\frac{1}{n}
\sum^n_\zs{j=[\sqrt{n}]+1}\,\xi^{2}_\zs{j}
=\sigma
\frac{n-[\sqrt{n}]}{n}
+ \frac{1}{\sqrt{n}}
\,
B_\zs{2,n}(x')
\,,
$$
where the function $B_\zs{2,n}(x')$ is defined in
\eqref{sec:Mo_Proof_B_M}
 for $x'_\zs{j}=1/\sqrt{n}\Chi_\zs{\{\sqrt{n}<j\le n\}}$.
Therefore, similarly to \eqref{sec:Mo.13_Ub}
we find
$$
\E_\zs{Q}
\left\vert
\frac{1}{n}
\sum^n_\zs{j=[\sqrt{n}]+1}\,\xi^{2}_\zs{j}
-\sigma
\right\vert
\le
\frac{\sigma\left(1+\sqrt{2+4\Pi(x^{4})}\right)}{\sqrt{n}}
\,.
$$
This implies that
\begin{equation*}\label{sec:Mo.15++-}
\E_\zs{Q}|\wh{\sigma}_\zs{n}-\sigma| \leq
\frac{4\|\dot{S}\|^2+2\sqrt{\sigma}\|\dot{S}\|
+
\sigma\left(1+\sqrt{2+4\Pi(x^4)}\right)}{\sqrt{n}}
\end{equation*}
and, therefore, we obtain the bound \eqref{sec:Mo.15++-}.
Hence Proposition \ref{sec:Mo.Prop.1}.
\endproof

\bigskip

\subsection{Proof of Proposition \ref{Pr.sec:Stc.1}}

\proof
Taking into account the definition
of $I_\zs{t}(f)$ in \eqref{sec:Imp.4} and  \eqref{sec:In.1+1}
we obtain through the Ito formula that
\begin{align}\label{sec:Stc.2}
I_\zs{t}(f)\,I_\zs{t}(g)=
\sigma\,
(f,g)_\zs{t}
+
\M^{f,g}_\zs{t}\,,
\end{align}
where
$$
\M^{f,g}_\zs{t}=\int^{t}_\zs{0}\,\Upsilon^{f,g}_\zs{s-}
\,
\d \xi_\zs{s}
+
\sigma^{2}_\zs{2}\,
\int^{t}_\zs{0}
f(s)\,g(s)\,
\d m_\zs{s}\,,
$$
$\Upsilon^{f,g}_\zs{s}=f(s) I_\zs{s}(g) + g(s) I_\zs{s}(f)$ and $m_\zs{t}=x^{2}*(\mu-\tilde{\mu})_\zs{t}$.
Using now the inequality \eqref{Novikov++} with $\Upsilon=x f$ and $p=2$ we obtain that
for any $f\in\L_\zs{2}[0,t]$
$$
\E\,\left( \int^{t}_\zs{0}f(s)\d z_\zs{s}\right)^{2}\le C\Vert f\Vert^{2}_\zs{t}<\infty
\,.
$$
Therefore, taking the expectation in
\eqref{sec:Stc.2} we obtain \eqref{sec:Stc.1}. Hence
Proposition \ref{Pr.sec:Stc.1}.
\endproof

\bigskip

\subsection{Property of Penalty term}

\begin{lemma}\label{Lem.A.1}
For any $n\geq 1$, $\lambda\in\Lambda$ and $0<\varepsilon<1$
\begin{equation}
\label{penalty-00}
P_\zs{n}(\lambda)\leq\frac{\E\,\Er_\zs{n}(\lambda)}{1-\varepsilon}+\frac{c^{*}_\zs{n}}{n\varepsilon(1-\varepsilon)}
\,,
\end{equation}
where $c^{*}_\zs{n}$ is defined in
\eqref{sec:Mo.9+1_psi}.
\end{lemma}
\proof
By the definition of $\Er_\zs{n}(\lambda)$ one has
\begin{align*}
\Er_\zs{n}(\lambda)&=\sum^{n}_\zs{j=1}\,(\lambda(j)\theta_\zs{j}^*-\theta_j)^2
=\sum^{n}_\zs{j=1}\,\left(\lambda(j)(\theta_\zs{j}^*-\theta_j)+(\lambda(j)-1)\theta_j\right)^2 \\[2mm]
&
\ge
\sum^{n}_\zs{j=1}\,\lambda(j)^2(\theta_\zs{j}^*-\theta_j)^2+
2\sum^{n}_\zs{j=1}\,\lambda(j)(\lambda(j)-1)\theta_j(\theta_\zs{j}^*-\theta_j).
\end{align*}
Taking into account the condition $\B_\zs{2}$) and the definition
\eqref{sec:Imp.12_Imp} we obtain that the last term in tho sum can be replaced as
$$
\sum^{n}_\zs{j=1}\,\lambda(j)(\lambda(j)-1)\theta_j(\theta_\zs{j}^*-\theta_j)
=
\sum^{n}_\zs{j=1}\,\lambda(j)(\lambda(j)-1)\theta_j(\wh{\theta}_\zs{j}-\theta_j)
\,,
$$
i.e.
$
\E\,\sum^{n}_\zs{j=1}\,\lambda(j)(\lambda(j)-1)\theta_j(\theta_\zs{j}^*-\theta_j)=0$ and, therefore, taking into account the definition
\eqref{sec:Mo.9_P} we obtain that
\begin{align*}
\E\,\Er_\zs{n}(\lambda)&\geq\sum^{n}_\zs{j=1}\,\lambda(j)^2\E\,(\theta_\zs{j}^*-\theta_j)^2=
\sum^{n}_\zs{j=1}\,\lambda(j)^2\E\,\left(\frac{\xi_\zs{j}}{\sqrt{n}}-g_\zs{\lambda}(j)\wh{\theta}_\zs{j}\right)^2\\[2mm]
&
\ge
P_\zs{n}(\lambda)-\frac{2}{\sqrt{n}}\E\,\sum^{n}_\zs{j=1}\,\lambda(j)^2g_\zs{\lambda}(j)\wh{\theta}_\zs{j}
\xi_\zs{j}
\\[2mm]&
\ge (1-\varepsilon)\,P_{n}(\lambda)
- \frac{1}{\varepsilon}\E\,\sum^{n}_\zs{j=1}\,g_\zs{\lambda}^{2}(j)\wh{\theta}_\zs{j}^{2}
\,.
\end{align*}
The inequality \eqref{sec:Mo.13_Ub++c-n}
 implies the bound \eqref{penalty-00}. Hence  Lemma \ref{Lem.A.1}.
\endproof

\subsection{The absolute continuity of distributions for the L\'evy processes}\label{subsec:App.5++}

In this section we study the absolute continuity for the  the L\'evy processes defined as

\begin{equation}\label{sec:App.5++.1}
  \d y_t=S(t)\d t+\d \xi_t\,,
  \quad 0\le t\le T\,,
 \end{equation}
where $S(\cdot)$ is any arbitrary nonrandom square integrated function, i.e. from $\L_\zs{2}[0,T ]$ and
$(\xi_\zs{t})_\zs{0\le t\le T}$ is the L\'evy process of the form \eqref{sec:In.1+1}
with nonzero constant $\sigma_\zs{1}$.
 We  denote by $\P_\zs{y}$ and $\P_\zs{\xi}$ the distributions of the processes
$(y_\zs{t})_\zs{0\le t\le 1}$
and
$(\xi_\zs{t})_\zs{0\le t\le 1}$
on the Skorokhod space $\D[0,T]$.  Now for any $(x_\zs{t})_\zs{0\le t\le T}$
from $\D[0,T]$
 we set
\begin{equation}\label{sec:App.++.1}
\Upsilon_\zs{T}(x)=
\exp\left\{\int^{T}_\zs{0}\,\frac{S(u)}{\sigma^{2}_\zs{1}}\,\d x^{c}_\zs{u}
-\,\int^{T}_\zs{0}\,
\frac{S^{2}(u)}{2\sigma^{2}_\zs{1}}\,
\d u
\right\}
\,,
\end{equation}
where $x^{c}$ is the continuous part of the process $x$ defined in \eqref{sec:App.7}.
Now we study the measures $\P_\zs{y}$ and
$\P_\zs{\xi}$ in $\D[0,T]$.
\begin{proposition}\label{Pr.sec:App.1++}
For any  $T>0$ the measure $\P_\zs{y}\ll \P_\zs{\xi}$
in $\D[0,T]$
and the Radon-Nikodym derivative
is
$$
\frac{\d\P_\zs{y}}{\d\P_\zs{\xi}}(\xi)
=\Upsilon_\zs{T}(\xi)
\,.
$$
\end{proposition}
\noindent {\bf  Proof.}
Note that to show  this proposition it suffices to check that
for any $0=t_\zs{0}<\ldots<t_\zs{n}=T$
any $b_\zs{j}\in\bbr$ for $1\le j\le n$
$$
\E\,\exp\left\{i\sum^{n}_\zs{l=1}b_\zs{j}(y_\zs{t_\zs{j}}-y_\zs{t_\zs{j-1}})\right\}
=
\E\,\exp\left\{i\sum^{n}_\zs{l=1}b_\zs{j}(\xi_\zs{t_\zs{j}}-\xi_\zs{t_\zs{j-1}})\right\}\Upsilon_\zs{T}(\xi)
\,.
$$
taking into account that the processes $(y_\zs{t})_\zs{0\le t\le T}$ and
$(\xi_\zs{t})_\zs{0\le t\le T}$
have the independent homogeneous increments, to this end one needs to check only
that for any $b\in\bbr$ and $0\le s<t\le T$
\begin{equation}\label{sec:App.++.2}
\E\,\exp\left\{i\,b (y_\zs{t}-y_\zs{s})\right\}
=
\E\,\exp\left\{i\,b (\xi_\zs{t}-\xi_\zs{s})\right\}\frac{\Upsilon_\zs{t}(\xi)}{\Upsilon_\zs{s}(\xi)}
\,.
\end{equation}
To check this equality note that
the process
$$
\Upsilon_\zs{t}(\xi)=\exp
\left\{\int^{t}_\zs{0}\,\frac{S(u)}{\sigma_\zs{1}}\,\d w_\zs{u}
-\,\int^{t}_\zs{0}\,
\frac{S^{2}(u)}{2\sigma^{2}_\zs{1}}\,
\d u
\right\}
$$
is the gaussian martingale. From here we directly obtain the squation \eqref{sec:App.++.2}. Hence
Proposition \ref{Pr.sec:App.1++}.
\endproof

\bigskip
\bigskip

\subsection{Properties of the term $\R_\zs{0,n}$}

\begin{lemma}\label{Le.sec:App.3+1}
For any $m>0$ the term $\R_\zs{0,n}$ introduced in \eqref{sec:Lo.12}
satisfies the following property
\begin{equation}\label{sec:App.5++R++0}
\lim_\zs{n\to\infty}\,n^{m}\,
\R_\zs{0,n}\,=0\,.
\end{equation}
\end{lemma}
\proof
First, setting $\zeta_\zs{n}=\sum^{d_\zs{n}}_\zs{j=1}\,\kappa^{2}_\zs{j}\,a_\zs{j}$,
we obtain that
$$
\left\{
S_\zs{\kappa}\notin W_\zs{k,r}
\right\}
=
\left\{
\sum^{d_\zs{n}}_\zs{j=1}\,\kappa^{2}_\zs{j}\sum^{k}_\zs{l=0}\,
\Vert \phi^{(l)}_\zs{j}\Vert^{2}
> r
\right\}
=
\left\{
\zeta_\zs{n}
> r
\right\}
\,.
$$
Moreover, note that one can check directly that
$$
\lim_\zs{n\to \infty}\,
\E\,\zeta_\zs{n}=
\lim_\zs{n\to \infty}\,
\frac{1}{v_\zs{n}}
\sum^{d_\zs{n}}_\zs{j=1}\,s^{*}_\zs{j}\,a_\zs{j}=r_\zs{\ve}=
(1-\ve)r\,.
$$
So, for sufficiently large $n$ we obtain that
$$
\left\{
S_\zs{\kappa}\notin W_\zs{k,r}
\right\}
\subset
\left\{
\wt{\zeta}_\zs{n}>
r_\zs{1}
\right\}
\,,
$$
where $r_\zs{1}=r\ve/2$,
$$
\wt{\zeta}_\zs{n}=\zeta_\zs{n}-\E\,\zeta_\zs{n}
=\frac{1}{v_\zs{n}}\,\sum^{d_\zs{n}}_\zs{j=1}\,s^{*}_\zs{j}a_\zs{j}\wt{\eta}_\zs{j}
\quad\mbox{and}\quad
\wt{\eta}_\zs{j}=\eta^{2}_\zs{j}-1\,.
$$
Through the correlation inequality from
\cite{GaPeSPA_2013}
we can get that for any $p\ge 2$ there exists some constant $C_\zs{p}>0$  for which
$$
\E\,\wt{\zeta}^{p}_\zs{n}\le C_\zs{p}
\frac{1}{v^{p}_\zs{n}}\,
\left( \sum^{d}_\zs{j=1}\,
(s^{*}_\zs{j})^{2}a^{2}_\zs{j}
\right)^{p/2}
\le C \,v^{-\frac{p}{4k+2}}_\zs{n}
\,,
$$
i.e.  the expectation
$\E\,\wt{\zeta}^{p}_\zs{n}\to 0$ as $n\to\infty$. Therefore, using the Chebychev inequality
we obtain that for any $m>1$
$$
n^{m}\P(\wt{\zeta}_\zs{n}> r_\zs{1})\to 0
\quad\mbox{as}\quad
n\to\infty\,.
$$

\noindent
Hence Lemma \ref{Le.sec:App.3+1}. \endproof


\bigskip

\bigskip



\begin{thebibliography}{100}


\bibitem{Akaike1974} Akaike H. A new look at the statistical model identification.
{\em IEEE Trans. on Automatic Control} \textbf{19} (1974) 716--723.

\bibitem{BarndorffNielsenShephard2001}
O. E. Barndorff-Nielsen and N. Shephard.
Non-Gaussian Ornstein-Uhlenbeck-based models and some of their uses in
financial mathematics. {\em J. Royal Stat. Soc.} \textbf{B 63} (2001) 167--241.


\bibitem{BarronBirgeMassart1999}
Barron  A.,  Birg\'e L. and  Massart P. (1999) Risk bounds for model selection via penalization.
{\em Probab. Theory Relat. Fields} \textbf{113},  301--415.


\bibitem{BichtelerJacod1983}
Bichteler K., Jacod J.  {\em Calcul de Malliavin pour les diffusions avec sauts: existence d’une densit\'e dans le cas unidimensionnel}.  S\'eminaire de probabilit\'e, XVII, Lecture Notes in Math., {\bf  986}, Springer, Berlin, 1983, 132--157.


\bibitem{Bertoin1996}
J. Bertoin. {\em L\'evy Processes.}
Cambridge University Press, Cambridge, 1996.











\bibitem{ComteGenenCatalot2011}
Comte, F. and Genon-Catalot, V. (2011) Estimation for L\'evy processes from high frequency data within a
  long time interval. {\em The Annals of Statistics}, {\bf 39}(2), 803 -- 837.





\bibitem{ContTankov2004}
Cont R., Tankov P. {\em Financial Modelling with Jump Processes.}
Chapman \& Hall, 2004.


\bibitem{FourdrinierStrawderman1996}
Fourdrinier D., Strawderman W. E. (1996). A paradox concerning shrinkage estimators: should a known scale parameter be replaced by an estimated value in the shrinkage factor? {\sl Journal of Multivariate Analysis}, 59(2), 109 --140.

\bibitem{FourdrinierPergamenshchikov2007}
Fourdrinier D., Pergamenshchikov S.  (2007)
Improved selection model method for the regression with dependent noise.
 {\em Annals of the Institute of Statistical Mathematics},
 \textbf{59}(3), 435--464.

\bibitem{Flaksman2002}
Flaksman, A.G. (2002) Adaptive signal processing in antenna arrays with allowance for the rank of the impule-response matrix of a multipath channel
  {\em Radiophysics and Quantum Electronics}, {\bf 45} (12),  977 -- 988.







\bibitem{GaltchoukPergamenshchikov2006}
Galtchouk L.I., Pergamenshchikov S.M. (2006)
Asymptotically efficient estimates for nonparametric regression
  models. {\em Statistics and Probability Letters}, {\bf 76} (8), 852--860.


\bibitem{GaltchoukPergamenshchikov2009a}
Galtchouk L.I., Pergamenshchikov S. M. (2009) Sharp non-asymptotic oracle inequalities
for nonparametric heteroscedastic regression models. {\em Journal of Nonparametric
Statistics} {\bf 21} (1),  1 - 16.



\bibitem{GaltchoukPergamenshchikov2009b}
Galtchouk L., Pergamenshchikov S. (2009) Adaptive
asymptotically efficient estimation in heteroscedastic
nonparametric regression.{\em Journal of Korean Statistical
Society}, \textbf{38}(4), 305--322.


\bibitem{GaPeSPA_2013}
Galtchouk L., Pergamenshchikov S. (2013)
Uniform concentration inequality for ergodic diffusion processes observed at discrete times.
{\em Stochastic Processes and their Applications}, {\bf 123},  91–109



\bibitem{IbragimovKhasminskii1981}
Ibragimov I. A., Khasminskii R. Z.
{\em Statistical Estimation: Asymptotic Theory}.
Springer,  New York, 1981.

\bibitem{JacodShiryaev2002} Jacod J., Shiryaev A.N. {\em Limit theorems for stochastic processes.} 2nd edition, Springer, Berlin, 2002.

\bibitem{JamesStein1961}
James W., Stein C. (1961). Estimation with quadratic loss. {\sl In Proceedings of the Fourth Berkeley Symposium Mathematics, Statistics and Probability, University of California Press, Berkeley}, 1, 361–380






\bibitem{Kassam1988}
Kassam S.A. Signal detection in non-Gaussian noise. – New York: Springer-Verlag Inc.,  IX,
1988.






\bibitem{KonevPergamenshchikov2009a}
Konev V. V., Pergamenshchikov S. M.  Nonparametric
estimation in a  semimartingale regression model. Part 1. Oracle
Inequalities. {\em
 Journal of Mathematics and Mechanics
of Tomsk State University} \textbf{3} (2009) 23--41.

\bibitem{KonevPergamenshchikov2009b}
Konev V. V., Pergamenshchikov S. M.
Nonparametric estimation in a  semimartingale regression
model. Part 2. Robust asymptotic efficiency.
{\em
 Journal of Mathematics and Mechanics
of Tomsk State University}
 \textbf{4} (2009) 31--45.

\bibitem{KonevPergamenshchikov2010}
Konev V. V., Pergamenshchikov S. M. General model
selection estimation of a periodic regression with a Gaussian
noise. {\sl Annals of the Institute of Statistical Mathematics} \textbf{62} (2010) 1083--1111.


\bibitem{KonevPergamenshchikov2012}
Konev V. V., Pergamenshchikov S. M.
 Efficient robust nonparametric estimation in a semimartingale
  regression model. {\em Ann. Inst. Henri Poincar\'e Probab. Stat.}, \textbf{48} (4), 2012, 1217--1244.


\bibitem{KonevPergamenshchikov2015}
Konev V. V., Pergamenshchikov S. M.
 Robust model selection for a semimartingale continuous time regression from
 discrete data. {\em Stochastic processes and their applications}, \textbf{125}, 2015, 294 -- 326.



\bibitem{KPP2014}
Konev V., Pergamenshchikov S. and Pchelintsev E. (2014)
Estimation of a regression with the pulse type noise from discrete data.
{\sl Theory Probab. Appl.}, {\bf 58} (3),  442--457.

\bibitem{Kutoyants1977}
Kutoyants Yu. A.
Estimation of the signal parameter in a Gaussian
Noise.
{\sl Problems of Information Transmission}, 1977,  vol. 13 (4), p. 29 -- 36.

\bibitem{Kutoyants1984}
Kutoyants Yu. A. {\em Parameter Estimation for Stochastic Processes}.
Heldeman-Verlag, Berlin, 1984.


\bibitem{JacodShiryaev1987}
Jacod J.,  Shiryaev A. N. {\em Limit theorems for stochastic processes}.
Vol.1, Springer, New York, 1987.

\bibitem{Mallows1973}
Mallows C. Some comments on $C_\zs{p}$.{\em Technometrics} \textbf{15} (1973) 661--675.



\bibitem{MarchallOlkin1979}
 Marshall A.W.,  Olkin I.
{\sl Inequalities: Theory of Majorization and Its Applications.}
Springer Series in Statistics,
 Springer,
 Academic.



\bibitem{MarinelliRockner2014}
Marinelli C., R\"ockner M. (2014)
On maximal inequalities for purely discontinuous martingales in infinite dimensions.
{\em S\'eminaire de Probabilit\'es}, Lect. Notes Math., {\bf XLVI} , 293--315.

\bibitem{Novikov1975}
Novikov A.A. (1975) On discontinuous martingales.  {\em Theory Probab. Appl.}, {\bf 20},
1, 11--26.



\bibitem{Nussbaum1985}
Nussbaum M. (1985) Spline smoothing in regression models and
asymptotic efficiency in $\L_2$.- {\em Ann. Statist.} \textbf{13},  984--997.

\bibitem{Pchelintsev2013}
Pchelintsev E. (2013)
Improved estimation in a non-Gaussian parametric regression.
{\sl Stat. Inference Stoch. Process.}, {\bf 16} (1),  15 -- 28.



\bibitem{Pinsker1981}
Pinsker M.S.  (1981) Optimal filtration of square integrable
signals in gaussian white noise. {\sl Problems of Transimission
information} \textbf{17}, 120--133.



\bibitem{Proakis1995}
Proakis J. G. {\em Digital Communications}.
McGraw-Hill, New York, 1995.

















\end{thebibliography}
\end{document}